\documentclass[11pt,a4paper]{article}
\usepackage[applemac]{inputenc}

\usepackage{amsfonts}
\def\up#1{\raise 1ex\hbox{\small #1}}
\usepackage[matrix,arrow]{xy}
\xyoption{all}

\begin{document}
\bigskip

\centerline{ \large \bf Sur la $K$-th\'eorie du foncteur norme}


\bigskip

{\centerline{ \large \bf Max Karoubi\footnote{Universit\'e Paris 7-Denis Diderot, UMR 7586 du CNRS, Case 7012,  175, rue du Chevaleret, 75205 Paris Cedex 13, France. Courriel : max.karoubi@gmail.com} \& Thierry Lambre\footnote{Universit\'e Blaise Pascal, UMR 6620 du CNRS, Les C\'ezeaux, 63177 Aubi\`ere, Cedex, France. Courriel : thierry.lambre@math.univ-bpclermont.fr }}}

\bigskip

{\bf R\'esum\'e :}

La $K$-th\'eorie d'un foncteur peut \^etre vue comme la version relative de celle d'un anneau unitaire. Aprs une description du groupe $K_0$ d'un foncteur, nous montrons que cette $K$-th\'eorie de foncteur poss\`ede  une suite exacte longue de Mayer-Vietoris.

Dans le cas d'une extension galoisienne de  corps de nombres  $F/L$, d'anneaux d'entiers respectifs $A$ et $B$, le groupe $K_0$ du ``foncteur norme'' est une extension d'un sous-groupe du   groupe des classes d'idaux $Cl(A)$ du corps $F$ par le groupe de cohomologie de Tate $\widehat H^0(G,A^*)$.

La suite exacte longue de Mayer-Vietoris permet d'expliciter un quotient du  sous-groupe 
$${}_NCl(A):=\ker N: Cl(A)\to Cl(B)$$ du groupe des classes $Cl(A)$, o $N$ dsigne la norme,   sous la forme d'une suite exacte courte 
$$\xymatrix{
 1\ar[r]& 
 B^*/B^*\cap N(F^*)\ar[r]^-{}&
 \widehat H^0(G, F^*)\cap \widehat H^0(G, {\cal U}_F)\ar[r]&
  {}_NCl(A)/I_GCl(A)\ar[r]& 1
 }$$
 o ${\cal U}_F$ est le groupe des units semi-locales du corps $F$. Pour conclure cet article, nous proposons quelques applications arithmtiques.

\medskip
{\bf Abstract:}

The $K$-theory of a functor may be viewed as a relative version of the $K$-theory of a ring. After its description, we prove a Mayer-Vietoris exact sequence in this framework. 

In the case of a Galois extension of a number field $F/L$ with rings of integers $A$, $B$ respectively, this $K$-theory of the ``norm functor'' is an extension of a subgroup of the ideal class group $Cl(A)$ of $F$ by the Tate cohomology group $\widehat H^0(G, A^*)$.

The Mayer-Vietoris exact sequence enables us to describe in a quite explicit way a quotient of the subgroup 
$${}_NCl(A):=\ker N: Cl(A)\to Cl(B)$$
of the ideal class group $Cl(A)$, where $N$ is the norm. We also prove  a short exact sequence
$$\xymatrix{
 1\ar[r]& 
 B^*/B^*\cap N(F^*)\ar[r]^-{}&
 \widehat H^0(G,F^*)\cap \widehat H^0(G,{\cal U}_F)\ar[r]&
  {}_NCl(A)/I_GCl(A)\ar[r]& 1
 }$$
 where ${\cal U}_F$ is the group of semi-local units of $F$.
 
Finally, we conclude this paper by  applications of our methods to Number Theory.

\medskip

Classification A.M.S.
Primaire : 19F05, 11R70,  11R34.\ 
Secondaire : 11R37, 11R29.

\bigskip

{\bf Introduction.}

\bigskip

Cet article a pour origine notre souhait d'utiliser dans un cadre arithmtique des outils de $K$-thorie en bas degr dont les ides remontent  H. Bass ([B]).

Les groupes $K_0$ et $K_1$ d'un anneau unitaire $A$ ou, ce qui revient au m\^eme, de la catgorie ${\bf Proj}(A)$,  ont t intensivement tudis ds leur apparition. Dans cet article nous tudions  un groupe moins tudi et fortement reli aux deux prcdents, le groupe de $K$-thorie d'un foncteur additif. Le groupe $K_0$ d'un tel foncteur $\varphi:{\bf Proj}(A)\to {\bf Proj}(B)$ appara\^{\i}t comme un invariant assez fin pour permettre  l'observation  de phnomnes que les groupes $K_0$ et $K_1$  eux seuls ne semblent pas dtecter.

\smallskip

Cet article dbute par une description minutieuse du groupe 
$K_0$ du foncteur norme ${\bf N}_{A/B}$ associ  une extension  de corps de nombres $F/L$ d'anneaux d'entiers $A={\cal O}_F$ et $B={\cal O}_L$. 

Soit  
${\cal I}_F$  le groupe des  idaux fractionnaires de $A$.
Soient 
$N:F\to L$ et 
$N:{\cal I}_F\to {\cal I}_L$ 
les normes respectives sur $F$ et ${\cal I}_F.$ Dfinissons
${\cal N}: F^*\to L^*\times {\cal I}_F$ par 
${\cal N}(z)=(N(z), zA)$.
 On pose 
 $${\cal K}_0({\bf N}_{A/B})=\{(t,I)\in L^*\times {\cal I}_F\ \mid \ tB=N(I)\}/\hbox{Im\ }{\cal N}.$$

 Pour $(t,I)$ dans $L\times {\cal I}_F$, on d\'esigne par  $[t,I]\in {\cal K}_0({\bf N}_{A/B})$ la classe de $(t,I)$ mod $\hbox{\ Im}({\cal N})$.
\medskip

{\bf Th\'eor\`eme 3.4.}
{\sl Soit $F/L$ une extension  de corps de nombres.
Le groupe 
$K_0({\bf N}_{A/B})$   du foncteur  norme ${\bf N}_{A/B}$ est isomorphe au groupe 
${\cal K}_0({\bf N}_{A/B})$.
De plus, il existe une suite exacte  
$$\xymatrix{
A^*\ar[r]^-{N_{F/L}}&B^*\ar[r]^-{\sigma}&
{\cal K}_0({\bf N}_{A/B})\ar[r]^-{\rho}&Cl(A)\ar[r]^-{N}&Cl(B)
}$$ 
o\`u les homomorphismes $\sigma$ et $\rho$ sont respectivement donn\'es par $\sigma(v)=[v,A]$ et $\rho([t,I])=[I]$.}
\medskip

Dans un second temps, nous construisons  une suite exacte longue de Mayer-Vietoris pour la $K$-thorie d'un tel foncteur.

  \medskip
  
 {\bf Thorme 4.4.}
{\sl Soit $F/L$ une extension de corps de nombres, d'anneaux d'entiers $A$ et $B$ et d'anneaux d'adles ${\bf A}_F$ et ${\bf A}_L$ respectivement. On dsigne  par $\widehat A$ l'anneau des localiss complts de $A$  toutes les places finies de $F$. 

Soit $\varphi:{\bf Proj}(A)\to  {\bf Proj}(B)$ un foncteur additif et soient  $\varphi_{F/L}$, $\varphi_{\widehat A/\widehat B}$ et $\varphi_{{\bf A}_F/{\bf A}_L}$ les foncteurs induits par  $\varphi$. Ces foncteurs s'insrent alors dans une suite exacte longue
$$\xymatrix{\cdots \ar[r]& 
K_{r+1}(\varphi_{{\bf A}_F/{\bf A}_L})\ar[r]^-{\partial}& 
K_r(\varphi)\ar[r]&
K_r(\varphi_{F/L})\times K_r(\varphi_{\widehat A/\widehat B})\ar[r]&
K_r(\varphi_{{\bf A}_F/{\bf A}_L})\ar[r]^-{\partial}& \cdots}$$
}

\medskip
 Dans le cas du foncteur norme et pour $r=0$, cette suite exacte peut \^etre prcise.
 Introduisons pour cela une notation : si $N:X\to Y$ est un morphisme de groupes, nous notons 
 $${}_NX:=\ker N$$ et
 $$Y_N:=\hbox{coker N}$$
 \medskip

{\bf Thorme 5.4.} (suite exacte des noyaux-conoyaux)
{\sl Soit $F/L$ une extension de corps de nombres d'anneaux d'entiers respectifs $A$ et $B$. Soient
 ${\cal U}_F$,  ${\cal U}_L $ les groupes d'units semi-locales de $F$ et $L$ respectivement,
${\bf J}_F$, ${\bf J}_L$ les groupes d'idles  de $F$ et $L$ respectivement. Notons indiffremment par la m\^eme lettre $N$ les diffrentes normes $A^*\to B^*$, 
$F^*\to L^*$, 
${\cal U}_F\to {\cal U}_L$ et ${\bf J}_F\to {\bf J}_L$.
\smallskip

On a alors une suite exacte  
$$\xymatrix{
&1\ar[r]&
{}_NA^*\ar[r]^-{\Delta}&
{}_NF^*\times{}_N{\cal U}_F\ar[r]^-{\mu_1}&
{}_N{\bf J}_F\ar[dldl]^-{\partial}&\cr
\cr
&&K_0({\bf N}_{A/B})\ar[r]^-{i}&
L^*_N\times ({\cal U}_L)_N\ar[r]^-{\mu}&
({{\bf J}_L})_N\ar[r]^-{\partial'}&
Cl(B)_N\ar[r]&1.}$$
} 

\medskip

 Cet article s'achve par quelques calculs explicites dans le cas des extensions cycliques.

\medskip

{\bf Thorme 6.6.}
{\sl Soit $F/L$  une extension  cyclique de corps de nombres d'anneaux d'entiers $A$ et  
 $B$ respectivement, de groupe de Galois $G$. On pose 
 $${}_NCl(A)=\ker N: Cl(A)\to Cl(B).$$
 
 Soit $I_GCl(A)$ le sous-groupe du groupe des classes $Cl(A)$ engendr par les lments de la forme $(1-g)[I]$, $g\in G$, $[I]\in Cl(A)$.  On a alors  un isomorphisme de groupes 

$${}_NCl(A)/I_G Cl(A)\ \cong
\ L^*/B^*N(F^*)\ \ \cap \ \ {\cal U}_L/B^*N({\cal U}_F).$$

De plus, on a une suite exacte 
 
 $$\xymatrix{
 1\ar[r]& 
 B^*/B^*\cap N(F^*)\ar[r]^-{}&
 \widehat H^0(G,F^*)\cap \widehat H^0(G,{\cal U}_F)\ar[r]&
  {}_NCl(A))/I_G Cl(A)\ar[r]& 1,
 }$$
o $\widehat H^0(G,-)$ dsigne le groupe $H^0$ de cohomologie de Tate.}
\medskip

Nous en dduisons 

\medskip

{\bf Thorme 6.8.}
{\sl Soit $F/L$ une extension cyclique de corps de nombres d'anneaux d'entiers $A={\cal O}_F$ et $B={\cal O}_L$. On suppose l'extension ramifie de degr premier $\ell$, de groupe de Galois $G$. Dsignons par $t$ le nombre de places finies de $L$ ramifies dans $F$ et par
$I_GCl(A)$ le sous-module de $Cl(A)$ engendr par les lments de la forme $(1-g)a$, $a\in Cl(A)$. Posons  
$${}_NCl(A):=\ker N:Cl(A)\to Cl(B).$$ Alors 
${}_NCl(A)/I_GCl(A)$ est un ${\bf F}_\ell$-espace vectoriel de dimension infrieure ou gale  $t-1$.}

\medskip

En particulier, pour une extension $F/{\bf Q}$ cyclique, si on dsigne  par $Cl(A)_G$ le groupe des co-invariants du groupe des classes de $A$ et par 
$\hat {\bf Z}^*={\cal U}_{\bf Q}=\prod_{p} {\bf Z}_p^*$,  on a un isomorphisme
  de groupes
$$Cl(A)_G\ \cong\  
 {\bf Q}^*/{\bf Z}^* N(F^*)\ \cap\  \widehat{\bf Z}^*/{\bf Z}^* N({\cal U}_F).$$
\medskip

 Cette suite exacte permet par exemple de proposer une dmonstration nouvelle du rsultat suivant sur le $\ell$-rang du groupe des classes d'une extension cyclique d'ordre $\ell$. Ce rsultat est dj connu mais de dmonstration toujours dlicate.
  
  \medskip

{\bf Th\'eor\`eme 6.9.} ([H],[L])
{\sl Soit $F/{\bf Q}$ une extension cyclique de degr\'e premier $\ell$, de groupe de Galois $G$. Notons $Cl(A)_G$ le groupe des co-invariants du groupe des classes $Cl(A)$ sous l'action de son groupe de Galois.

Soit $t$ le nombre de diviseurs premiers du discriminant du corps $F$.
Alors le $\ell$-rang du groupe des co-invariants du groupe des classes du corps $F$ est \'egal \`a $t-1$, sauf si $\ell=2$, $F$ quadratique r\'eel et s'il existe un nombre premier $p$ congru \`a $3$ modulo $4$ et  divisant le discriminant du corps $F$. 

Dans ce dernier cas,  le $2$-rang du groupe des co-invariants du groupe des classes du corps $F$ est \'egal \`a $t-2$.}

\bigskip

Nous tenons  remercier chaleureusement T. Nguyen Quang Do (Universit de Franche-Comt) pour sa lecture attentive d'une version prliminaire de ce texte. Nous  remerions galement J. Berrick et Meng Fai Lim (Universit de Singapore)  pour nous avoir communiqu la remarque 5.6 ainsi que C. Movahhedi (Universit de Limoges) pour une correspondance fructueuse concernant le rsultat 6.8.
\bigskip

Ce texte est organis\'e comme suit

\smallskip

1. Notations et rappels.

\smallskip

2. Le groupe $K_0$  d'un foncteur.
\smallskip

3. Le foncteur norme 
\smallskip

4. La suite exacte longue de Mayer-Vietoris de la $K$-thorie d'un foncteur.
\smallskip

5. Le cas des corps de nombres et du foncteur norme.

\smallskip

6. Quelques applications.
\bigskip

{\bf 1. Notations et rappels.}
\medskip

Dans ce texte, nous considrons une extension $F/L$ de corps de nombres. Nous dsignons par :

$A={\cal O}_F$ et $B={\cal O}_L$  les anneaux d'entiers respectifs de $F$ et $L$, 

${\cal I}_F$  le groupe des  idaux fractionnaires de $A$,

$Cl(A)$ le groupe des classes d'idaux de $A$ (ou $F$),

$V(F)$  l'ensemble des places {\sl finies} de $F$,

 Soit $w$  une place finie de $F$ et soit $v$  une place finie de $L$. On note
$A_w$ (resp. $F_w$)   le complt de l'anneau $A$ (resp. du corps $F$)  la place finie $w$ de $F$. 

On dsigne de m\^eme par 

$\widehat A=\prod_{w\in V(F)}A_w$ et $\widehat B=\prod_{v\in V(L)}B_v$,

${\bf A}_F=\widehat A\otimes_A F$  l'anneau des adles de $F$, restreintes aux places finies,

${\cal U}_F=\widehat A^*=\prod_{w\in V(F)}A^*_w$   le groupe multiplicatif des units locales, c'est--dire les lments invresibles  de l'anneau $\widehat A$,

${\bf J}_F=\left({\bf A}_F\right)^*$   le groupe multiplicatif des idles, units de l'anneau ${\bf A}_F$.
\medskip

Pour un $G$-module $M$, on dsigne par $I_GM$ le sous-module de $M$ engendr par les lments de la forme $(1-g)m$, o $g$ appartient  $G$ et $m$ appartient  $M$. Le groupe $M_G=M/I_GM$ est le groupe des co-invariants du $G$-module $M$.
\medskip
 
Lorsque l'extension $F/L$ est galoisienne de groupe $G$, le groupe de cohomologie de Tate  $\widehat H^0(G,M)$ d'un $G$-module $M$ est dfini par
$$\widehat H^0(G, M)=M^G/\hbox{im\ } \nu,$$
o $\nu:M\to M$ est le morphisme de $G$-modules dfini par 
$\nu(m)=\sum_{g\in G}gm$.

Dans le cas galoisien, si $N:A^*\to B^*$ dsigne la norme,  on a donc 
$$\widehat H^0(G,A^*)=\hbox{coker\ }\left( N:A^*\to B^*\right)=B^*_N.$$ 

\medskip

Nous aurons besoin du  calcul de la $K$-thorie des anneaux $\widehat A$ et ${\bf A}_F$.
\medskip

{\bf Lemme 1.1.}
{\sl Soit $A$ l'anneau des entiers d'un corps de nombres $F$. On a alors des isomorphismes de groupes
$K_1(\widehat A)\cong{\cal U}_F$ et $K_1({\bf A}_F)\cong{\bf J}_F$.
}
\medskip
 
 {\sl Dmonstration.}
 L'anneau $\widehat A=\prod_{w\in V(F)}A_w$ est produit d'anneaux de valuations discrtes. Tous les anneaux $A_w$ sont de rang stable $d\leq 1$. De plus, pour tout $n\geq 3$, il existe $N$ tel que pour chaque place $w$ de $F$, tout  lment de $E_n(A_w)$ est produit d'au plus $N$ matrices lmentaires. Montrons que l'homomorphisme vident 
 $\phi_1: K_1(\widehat A)\to\prod_{w\in V(F)} K_1(A_w)$ est un isomorphisme.
 
 L'homomorphisme $\phi_1$ est injectif : soit $\alpha$ un lment de $GL(\widehat A)$ tel que pour chaque place $w$ de $F$, l'image $\alpha_w$ de $\alpha$ dans $K_1(A_w)$ soit gal  $1$. Pour $n$  assez grand, ceci implique que $\alpha_w$ est un lment de $E_n(A_w)$ et donc s'crit comme un produit de commutateurs 
 $\prod_{1\leq i\leq N}[\beta_{w,i},\  \gamma_{w,i}]$. D'aprs ce qui prcde, le nombre $N$ de commutateurs peut tre born  indpendamment de la place $w$.
 
  Pour $1\leq i\leq N$ introduisons les lments $b_i$ et  et $c_i$ de $\widehat A$ dfinis  par
 $$b_i=(\beta_{w,i})_{w\in V(F)} \hbox{\ \ et\ \ } 
 c_i=(\gamma_{w,i})_{w\in V(F)}.$$
 Dans $GL(\widehat A)$, on a l'galit
 $$\alpha=[b_1,c_1]\times \ \cdots \ \times [b_N,c_N],$$ ce qui signifie que la classe de $\alpha$ dans $K_1(\widehat A)$ est $1$.

 L'homomorphisme $\phi_1$ est surjectif : Puisque le rang stable des anneaux $A_w$ est born par $1$, il existe un entier $n$, indpendant de $w$  tel que tout lment de $K_1(A_w)$ puisse s'crire comme la classe d'une matrice $\alpha_w\in Gl_n(A_w)$. La famille des matrices $(\alpha_w)_{w\in V(F)}$ dfinit un lment de 
 $GL_n(\widehat A)$, et par suite un lment $\alpha$  de $K_1(\widehat A)$ dont l'image par $\phi_1$ est la classe  dans chaque $K_1(A_w)$  de la matrice $\alpha_w$. 
 
 Ce raisonnement  tablit un isomorphisme de groupes 
 $K_1(\widehat A)\cong\prod_{w\in V(F)}K_1(A_w)$. Puisque  $\prod_{w\in V(F)}K_1(A_w)={\cal U}_F$, le premier isomorphisme est tabli.

 \smallskip
 
 Cette dmonstration s'adapte sans difficult au cas du produit restreint ${\bf A}_F$. Pour cela, soit $S$ un ensemble fini de places de $F$. On pose traditionnellement 
$${\bf A}_F^S=\prod_{w\in S}F_w\times\prod_{w\not\in S} A_w$$ et 
$${\bf J}_F^S=\prod_{w\in S}F^*_w\times\prod_{w\not\in S} A^*_w$$
de sorte que ${\bf A}_F=\displaystyle\lim_{\to}{\bf A}_F^S$
et ${\bf J}_F=\displaystyle\lim_{\to}{\bf J}_F^S$ ({\sl cf.\ }[N]).

D'aprs les considrations prcdentes,  on a 
$$K_1({\bf A}_F^S)=\prod_{w\in S}K_1(F_w)\times\prod_{w\not\in S}K_1(A_w),$$
ce qui donne  
$K_1({\bf A}_F^S)={\bf J}_F^S$.

On sait que le foncteur $K_1$ commute aux limites inductives ({\sl cf.}\ par exemple [W2]) d'o
$$K_1({\bf A}_F)=\displaystyle\lim_\to K_1({\bf A}_F^S),$$
ce qui dmontre $K_1({\bf A}_F)={\bf J}_F$.

 \medskip
 {\bf Dfinition 1.2.}
 {\sl Soit $\Lambda$ un anneau. On dit qu'un lment $u$ de $K_0(\Lambda)$ est born par $N$ si
 $u=[E]-[\Lambda]^N$, o $E$ est un facteur direct de $\Lambda^{2N}$.
 \smallskip
 
Si $(\Lambda_i)_{i\in I}$ est une famille d'anneaux, 
le produit born des groupes $K_0(\Lambda_i)$, $i\in I$,
est le groupe not 
$$\prod_{i\in I}\!{}^{(b)}(K_0(\Lambda_i))$$
dont les lments sont les familles $u=(u_i)_{i\in I}$
avec   $u_i\in K_0(\Lambda_i)$ et telles qu'il existe
un entier $N=N(u)$ pour lequel $u_i$ est born par $N$ pour tout $i$.}
 \medskip

{\bf Proposition 1.3.}
{\sl Soit $A$ l'anneau des entiers d'un corps de nombres. Le morphisme de groupes  vident 
$$\phi_0: K_0(\widehat A)\to \prod_{w\in V(F)}K_0(A_w)$$
est injectif, d'image le produit born $\prod\!{}^{(b)}_{w\in V(F)}K_0(A_w)$. De plus, les groupes $K_0(\widehat A)$ et $K_0({\bf A}_F)$ sont isomorphes.}
\medskip

{\sl Dmonstration}. L'homomorphisme $\phi_0$ est d'image le produit born $\prod\!{}^{(b)}_{w\in V(F)}K_0(A_w)$. En effet,  puisque le rang stable de chaque $A_w$ est infrieur ou gal  $1$, tout module projectif  de type fini sur $A_w$ s'crit $E\oplus A_w^m$ o $E$ est de rang au plus $1$. Choisissons $N>1$. Soit 
$u=(u_w)_{w_n\in V(F)}$ un lment du produit born
$\prod_{w\in V(F)}^{(b)}K_0(A_w)$. Puisque chaque lment $u_w$ du produit born s'crit sous la forme $[E_w]-[A_w^N]$, o $E_w$ est facteur direct de $A_w^{2N}$, le module $E_w$ est l'image d'un projecteur $p_w$ de $\hbox{ Mat}_{2N}(A_w)$. Cette famille de projecteurs dfinit un projecteur $p$ de $Mat_{2N}(\widehat A)$, ce qui fournit un lment de $K_0(\widehat A)$ d'image $u$.

Montrons par ailleurs l'injectivit de $\phi_0$. Soit 
$\alpha=[E]-[\widehat A^N]$ un lment de $K_0(\widehat A)$ d'image nulle par $\phi_0$. On peut supposer que $E$ est l'image d'un projecteur $p$ de 
$\hbox{\ Mat}_{2N}(\widehat A)$. En remplaant le projecteur $p$ par $p\oplus 0$, ceci implique que pour toute place $w$ de $F$, $p_w\oplus 0$ est conjugu au projecteur $\hbox{\ diag} (1, 0)$ dans $\hbox{\ Mat}_{4N}(A_w)$ par une matrice inversible $\alpha_w$ de 
$\hbox{\ GL}_{4N}(A_w)$. On en dduit que $\alpha $ est nul dans $K_0(\widehat A)$.
\smallskip

Pour le calcul de $K_0({\bf A}_F)$, on remarque que si $S$ est un ensemble fini de places de $F$, le groupe 
$K_0(\prod_{w\not\in S}A_w)$ est isomorphe au produit born
$\prod_{w\not\in S}^{(b)}K_0(A_w)$. On en dduit 
$$K_0({\bf A}_F^S)\cong
\prod_{w\in S}K_0(F_w) \times
\prod_{w\not\in S}\!\!{}^{(b)}\ K_0(A_w).$$
Par passage  la limite inductive, on obtient 
$$K_0({\bf A}_F)\cong\prod_{w\in V(F)}\!\!\!\!\!\!{}^{(b)} \ K_0(A_w)\cong K_0(\widehat A).$$

\bigskip

{\bf 2. Le groupe $K_0$ d'un foncteur.}

\medskip

Soient $A$ et $B$ deux anneaux unitaires, non ncessairement commutatifs.   
Soit $\varphi:{\bf Proj}(A)\to {\bf Proj}(B)$ un foncteur additif et soit
$BQ(\varphi):BQ{\bf Proj}(A)\to BQ{\bf Proj}(B)$ l'application continue d\'efinie entre les espaces classifiants de la $K$-th\'eorie de Quillen des anneaux $A$ et $B$ [Q]. Soit ${\cal F}_\varphi$ la fibre homotopique de l'application $BQ(\varphi)$. 
\medskip

{\bf D\'efinition 2.1.}
{\sl Pour $r\geq 0$, le groupe  $K_r(\varphi):=\pi_{r+1}({\cal F}_\varphi)$ s'appelle le $r$-i\`eme groupe de $K$-th\'eorie du foncteur $\varphi$.}
\medskip

De cette d\'efinition on d\'eduit  la suite exacte longue tautologique

$$\xymatrix{
\cdots \ar[r]^-{\partial}& 
K_1(\varphi)\ar[r]&
K_1(A)\ar[r]^-{\phi_1}&
K_1(B)\ar[r]^-{\partial}&
K_0(\varphi)\ar[r]&
K_0(A)\ar[r]^-{\phi_0}&
K_0(B)&
}$$

avec $\phi_i=\pi_{i+1}BQ(\varphi).$

\medskip

{\sl  Description du groupe $K_0$  d'un foncteur.}
\medskip

Dans ce cadre de gnralit, le calcul de 
$K_r(\varphi)=\pi_{r+1}({\cal F}_\varphi)$ reste dlicat, y compris pour $r=0$. Voici une autre description de $K_0(\varphi)$.

Pour cela, rappelons une construction g\'en\'erale due \`a H. Bass ([B], VII.5). Soient ${\cal C}$ et ${\cal D}$ deux cat\'egories mono\"{\i}dales sym\'etriques et soit $\varphi:{\cal C}\to{\cal D}$ un foncteur mono\"{\i}dal et cofinal\footnote{L'exemple le plus important est celui d'une catgorie additive et d'un foncteur additif et cofinal.}. Le c\^one de $\varphi$ est la cat\'egorie $co(\varphi)$ constitu\'ee des triplets $(C,\alpha,C')$ o\`u $C$ et $C'$ sont des objets de ${\cal C}$ et o\`u $\alpha:\varphi (C)\to\varphi (C')$ est un isomorphisme dans la cat\'egorie ${\cal D}$. Un morphisme de $(C,\alpha,C')$ vers $(C_1,\alpha_1,C'_1)$ est un couple $(f,f')$ o $f:C\to C_1$ et $f':C'\to C'_1$ sont des morphismes dans ${\cal C}$
tels que 
$\alpha_1\circ \varphi(f)=\varphi(f')\circ\alpha$. L'ensemble des classes d'isomorphie d'objets de $co(\varphi)$ est un mono\"{\i}de ab\'elien dont on note $K(co(\varphi))$ le groupe de Grothendieck.

\medskip
{\bf D\'efinition 2.2.}
{\sl Le groupe $K_0(\varphi)$ de K-th\'eorie du foncteur $\varphi$ est le  quotient du groupe $K(co(\varphi))$ par le sous-groupe $R$ engendr\'e par les \'el\'ements de la forme 
$$(C,\alpha,C')+(C',\beta,C'')-(C,\beta\alpha,C'').$$

\smallskip

La classe de $(C,\alpha,C')$ modulo $R$ est not\'ee $[C,\alpha,C']$}.

\medskip
Le groupe $K_0(\varphi)$ d'un foncteur mono\"{\i}dal et cofinal s'ins\`ere dans 
la suite exacte de Bass qui s'\'ecrit 
$$\xymatrix{
K_1({\cal C})\ar[r]^-{\phi_1}
&K_1({\cal D})\ar[r]^-{\sigma}&K_0(\varphi)\ar[r]^-{\rho}&
K_0({\cal C})\ar[r]^-{\phi_0}& K_0({\cal D})}.\leqno (2.2)$$

\medskip

De plus, on sait ([Gr] et  [K3], p. 269) que si le foncteur $\varphi$ est additif et cofinal, le groupe $K_0(\varphi)$ d\'efini en 2.1 par $K_0(\varphi)=\pi_{1}({\cal F}_\varphi)$  
est isomorphe au  groupe d\'efini ci-dessus par g\'en\'erateurs et relations. 

\medskip

{\sl Description du groupe $K_0$ d'un foncteur entre catgories de Cartier.}
\medskip

L'exemple suivant nous sera trs utile dans la suite.
Soit $A$ un anneau commutatif unitaire intgre de corps des fractions $F$.
Soient $\hbox{Cart}(A)$ le groupe des idaux fractionnaires inversibles de l'anneau $A$ et $\hbox{div}_F:F^*\to \hbox{Cart}(A)$ le morphisme de groupes dfini par 
$\hbox{div}_F(z)=zA$.

Introduisons la catgorie $\hbox{{\bf Cart}}(A)$ dont les objets sont les idaux fractionnaires inversibles de $A$ et dont les morphismes sont donns par 
$\hbox{Hom}_{\hbox{{\bf Cart}}(A)}(I,J)=\hbox{Hom}_A(I,J)$;  la structure mono\"{\i}dale  est donne par le produit des idaux. Les objets de cette catgorie sont donc les lments du groupe de Cartier de $A$.

Soit $\varphi:{\bf Cart}(A)\to{\bf Cart}(B)$ un foncteur mono\"{\i}dal et cofinal et soient $F$ et $L$ respectivement les corps de fractions de $A$ et $B$. Considrons
le groupe produit fibr 
$L^*\times_\varphi\hbox{Cart}(A)$, c'est--dire le sous-groupe de $L^*\times\hbox{Cart}(A)$ constitu des lments $(t,I)$ tels que 
$tB=\varphi(I)$. 
Soit $\Phi:F^*\to L^*\times_\varphi\hbox{Cart}(A)$ le morphisme de groupes dfini par
$\Phi(z)=(\varphi(z),zA)$.

\medskip

{\bf Dfinition 2.3.}
\medskip

 {\sl On pose 
 $${\cal K}_0(\varphi)=
 \hbox{coker\ }(\Phi)=
 L^*\times_\varphi\hbox{Cart}(A)/\hbox{im\ }(\Phi).$$
\smallskip

Pour $(t,I)$ dans $L\times_\varphi\hbox{Cart}(A)$, on pose 
$[t,I]=(t,I)\ \hbox{mod\ }\hbox{im\ }\Phi$}.

\medskip

{ N. B.} On notera l'analogie entre  
${\cal K}_0(\varphi)$ et le groupe ${\cal U}(A;{\bf Z}/n)$ de [K-L], Thorme 1.1.

\medskip

Le groupe $K_0({\varphi})$ est constitu\'e d'\'el\'ements de la forme $[J,\alpha,I]$, o\`u $J$ et $I$ sont des id\'eaux fractionnaires de $F$ et o\`u $\alpha:\varphi(I)\to \varphi(I)$ est un isomorphisme de $B$-modules. En multipliant   $[J,\beta,I]$  par 
$1=[J^{-1},\ id_{\varphi(J^{-1})}, J^{-1}]$, on peut toujours \'ecrire un \'el\'ement de $K_0(\varphi)$ sous la forme $[A,\alpha, I]$. Par $B$-lin\'earit\'e, l'application $\alpha:B\to \varphi(I)$ est d\'etermin\'ee par $t=\alpha(1)\in L^*$. De plus, puisque $\alpha$ est un isomorphisme, on a $\varphi(I)=tB$.
\medskip

{\bf Lemme 2.4.}

{\sl 
L'application  qui \`a $[A,\alpha,I]$ de $K_0(\varphi)$ associe $ [\alpha(1),I]$ dans ${\cal K}_0(\varphi)$ est un isomorphisme.}
\medskip

{\sl Dmonstration}. 
L'application rciproque est celle qui  l'lment $[t,I]$ de 
${\cal K}_0(\varphi)$ associe l'lment $[A,\alpha, I]$ de $K_0(\varphi)$, o $\alpha:B\to \varphi(I)$ 
est l'isomorphisme de $B$-modules dfini par $\alpha(1)=t$.

\medskip

{\sl  Le  foncteur restriction des scalaires.}

\medskip

Soient $A$ et $B$ deux anneaux unitaires (non n\'ecessairement commutatifs). On suppose que $A$ est un $B$-module projectif de type fini. Notons 
${\bf res}:{\bf Proj}(A)\to {\bf Proj}(B)$ le foncteur restriction des scalaires  d\'efini pour tout $A$-module $P={}_AP$, de classe d'isomorphie $[{}_AP]$ par la formule ${\bf res}([{}_AP])=[{}_BP]$, o\`u  ${}_BP$ d\'esigne le $B$-module $P$ obtenu par restriction des scalaires.

Le foncteur ${\bf res}$ est \'evidemment additif. Puisque $A$ est un $B$-module projectif de type fini, il est \'egalement cofinal. En effet, pour un $B$-module projectif de type fini $Q$, on pose $P=A\otimes_BQ$ et on voit que ${}_BP$ est alors facteur direct d'un $Q^n$. Donc ${\bf res }([P])$ contient $[Q]$ en facteur direct.

\medskip
Supposons de plus que les anneaux $A$ et $B$. Les homomorphismes de la suite exacte de Bass
$$\xymatrix{
K_1(A)\ar[r]^-{res_1}&K_1(B)\ar[r]^-{\sigma}&
K_0({\bf res})\ar[r]^-{\rho}&
K_0(A)\ar[r]^-{res_0}&K_0(B)}$$

 sont donns respectivement par
$$res_1([{}_AP,\alpha])=[{}_BP,\alpha],$$ 
$$\sigma([Q,\beta])=[A\otimes_BQ,\ id_A\otimes_B\beta,\ A\otimes_BQ],$$  
$$\rho ([P,\alpha,P'])=[P]-[P']$$

et $$res_0([{}_AP])=[{}_BP].$$

\medskip

{\sl Description du morphisme de groupes $\hbox{res}_1:K_1(A)\to K_1(B)$.}

Pour dcrire $res_1$, rappelons une  construction du groupe  $K_1(A)$. Soit  ${\cal K}_1(A)$ le groupe ab\'elien libre engendr\'e par les classes d'isomorphie de couples $(P,\alpha)$ o\`u $P$ est un $A$-module projectif de type fini et o\`u $\alpha$ appartient \`a ${\bf Aut}_A(P)$ et soit ${\cal H}(A)$ le sous-groupe de  ${\cal K}_1(A)$ engendr\'e par les deux types d'\'el\'ements suivants:
$$(P,\alpha)+(P,\beta)-(P,\beta\alpha)$$
$$(P,\alpha)+(Q,\beta)-(P\oplus Q,\alpha\oplus\beta).$$ Le groupe  quotient ${\cal K}_1(A)/{\cal H}(A)$ est alors isomorphe   $K_1(A)$.

\medskip

{\bf Proposition 2.5.}
{\sl Soit ${\hbox{\sl d\'et}}_B:K_1(B)\to B^*$ le d\'eterminant.  
Pour $u\in A^*$, on a 
$${\hbox{\sl d\'et}}_B\circ \hbox{res}_1(u)= {\hbox{\sl d\'et}}(\mu_u)\in B^*$$ 
o\`u $\mu_u:A\to A$ est l'application $B$-lin\'eaire $\mu_u(z)=uz$.
\medskip

En particulier, si $A$ et $B$ sont des anneaux d'entiers de corps de nombres, d'aprs le thorme de Bass, Milnor et Serre [B-M-S], on a  $K_1(A)=A^*$ et $K_1(B)=B^*$. On en dduit  
 $\hbox{res}_1(u)= \hbox{dt}(\mu_u)$ pour tout $u\in K_1(A)$. 
}
\medskip

{\sl Dmonstration}. 
Pour $(P,\alpha)$ dans ${\cal K}_1(A)$ on note $[P,\alpha]=(P,\alpha)$ mod ${\cal H}(A)$. Avec ces notations, on a $$\hbox{res}_1([P,\alpha])=[{}_BP,\alpha].$$ 

Pour $u\in A^*$ donnons l'expression  de 
${\hbox{\sl d\'et}}_B\circ res_1(u)$ dans $K_1(B)$. 
On \'ecrit $u=[A,\mu_u]$ o\`u $\mu_u:A\to A$ est la multiplication par $u$. On a $res_1(u)=[{}_BA,\mu]$. Or $A$ est un $B$-module projectif de type fini. Si $A$ est facteur direct d'un $B$-module libre de rang $n$, on a donc
$${\hbox{\sl d\'et}}_B([{}_BA,\mu])=[\Lambda^n_B(B^n),\Lambda^n_B(\mu)]=[B,\gamma],$$
 o\`u $\gamma:B\to B$ est la multiplication par ${\hbox{\sl d\'et}}(\mu_u)$.
 
\bigskip

{\bf 3. Le  foncteur norme.}
\medskip

Supposons que $F/L$ soit une extension  finie de corps de nombres, d'anneaux d'entiers $A={\cal O}_F$, $B={\cal O}_L$. Soient $Cl(A)$ et $Cl(B)$ les groupes des classes d'idaux de $A$ et $B$ respectivement.

\medskip
{\bf D\'efinition 3.1.} 
{\sl Le foncteur restriction des scalaires {\bf res} est not 
$${\bf N}_{A/B}:{\bf Proj}(A)\to{\bf Proj}(B).$$ Il est d\'efini par  
 ${\bf N}_{A/B}([{}_AP])=[{}_BP]$.}

\medskip

Cette terminologie trouve sa justification dans le fait suivant :
sur le groupe $K_1$, le foncteur norme co\"{\i}ncide avec la norme $N_{F/L}$. En effet,  de la proposition 2.5, on dduit :
\medskip

{\bf Proposition 3.2.}
{\sl  L'homomorphisme $N_1:K_1(A)\to K_1(B)$ est donn\'e par 
$N_1(u)=N_{F/L}(u)$.}
\medskip

La suite exacte de Bass 2.2 du foncteur norme ${\bf N}_{A/B}:{\bf Proj}(A)\to{\bf Proj}(B)$ se r\'eduit  donc \`a 
 la suite exacte courte
$$\xymatrix{1\ar[r]&B^*_N\ar[r]^-{\sigma}&
K_0({\bf N}_{A/B})\ar[r]^-{\rho}&{}_NCl(A)\ar[r]&1.
 }\leqno (3.2)$$
o\`u  
$B^*_N=\hbox{coker\ }N:A^*\to B^*$ et
${}_NCl(A)=\ker N:Cl(A)\to Cl(B)$.
\medskip

Pour obtenir une description plus pr\'ecise du groupe $K_0({\bf N}_{A/B})$ au sein de cette extension, on utilise 2.4.
Introduisons pour cela
 le produit fibr\'e ${\cal K}=L^*\times_{{\cal I}_L} {\cal I}_F$,
o ${\cal I}_F$ est le groupe des idaux fractionnaires de $F$.
 
Le groupe ${\cal K}$ est donc le sous-groupe de  $L^*\times {\cal I}_F$ constitu\'e des \'el\'ements $(t,I)$ tels que l'id\'eal $N(I)$ soit principal de g\'en\'erateur $t$, c'est-\`a-dire $(t)=(N(I))$. Soit ${\cal N}:F^*\to {\cal K}$ le morphisme de groupes d\'efini par 
 ${\cal N}(z)=(N(z), zA)$. 
\medskip

{\bf D\'efinition 3.3.}
{\sl Posons $${\cal K}_0({\bf N}_{A/B})={\cal K}/Im({\cal N})=\{(t,I)\in L^*\times {\cal I}_F\ \mid \ tB=N(I)\}/\{(N(z),zA),z\in F^*\}.$$}
 \medskip
  
 Pour $(t,I)$ dans $\cal K$, d\'esignons par  $[t,I]\in {\cal K}_0({\bf N}_{A/B})$ la classe de $(t,I)$ mod $Im({\cal N})$.
\medskip

{\bf Th\'eor\`eme 3.4.}
{\sl Soit $F/L$ une extension  de degr\'e fini de corps de nombres.
Le groupe 
$K_0({\bf N}_{A/B})$ de $K$-th\'eorie  du foncteur  norme ${\bf N}_{A/B}$ est isomorphe au groupe 
${\cal K}_0({\bf N}_{A/B})$.
De plus, la suite exacte de Bass du foncteur  ${\bf N}_{A/B}$ s'identifie  la suite \hbox{exacte}  
$$\xymatrix{
A^*\ar[r]^-{N_{F/L}}&B^*\ar[r]^-{\sigma}&
{\cal K}_0({\bf N}_{A/B})\ar[r]^-{\rho}&Cl(A)\ar[r]^-{N}&Cl(B)
}$$ 
o\`u les homomorphismes $\sigma$ et $\rho$ sont donn\'es par $\sigma(v)=[v,A]$ et $\rho([t,I])=[I]$.}
\medskip

{\sl Dmonstration}. Considrons le diagramme suivant de catgories et de foncteurs, commutatif  isomorphisme prs.
$$\xymatrix{
{\bf Cart}(A)\ar[r]^-{{\bf N}_{Cart}}\ar[d]_-{{\bf i}_A}&{\bf Cart}(B)\ar[d]^-{{\bf i}_B}\cr
{\bf Pic}(A)\ar[r]^-{{\bf \Lambda}}&{\bf Pic}(B)\cr
{\bf Proj}(A)\ar[u]^-{{\hbox{\bf d\'et}}_A}\ar[r]^-{{\bf N}_{A/B}}&\ar[u]_-{{\hbox{\bf d\'et}}_B}{\bf Proj}(B)
}$$

La cat\'egorie ${\bf Pic}(A))$ a pour objets les $A$-modules projectifs de rang $1$. Le foncteur ${\bf \Lambda}$ est donn\'e par 
${\bf \Lambda}({}_AP)=\Lambda^{n}_{B}({}_BP))$, o $n=[F:L]$. 
Ce foncteur est cofinal (dans un contexte  multiplicatif). Par ailleurs, dans l'anneau de Dedekind $A$, si $P$ et $Q$ sont des modules projectifs de rang $1$, alors les modules $P\oplus Q$ et  $(P\otimes_A Q)\ \oplus\  A$ sont isomorphes ([M]). 
Or 
$\Lambda^{n}_ {B}({}_B(P\otimes_AQ\oplus A))=\Lambda^{n}_{B}({}_B(P\otimes_AQ))$
 tandis que 
$\Lambda^{n}_{B}({}_B(P\oplus Q))=\Lambda^{n}_{B}({}_BP)\otimes_B\Lambda^{n}_{B}({}_BQ)$. On a donc la relation 
$${\bf \Lambda}(P\otimes_AQ)={\bf \Lambda}(P)\otimes_B{\bf \Lambda}(Q),$$
ce qui montre que ce foncteur est additif au sens mono\"{\i}dal.
\medskip

On pose ${\hbox{\bf d\'et}}_A([P])=\Lambda^r_A(P)$ o\`u $r$ est le rang de $P$ en tant que $A$ module. Le thorme 3.4 r\'esulte alors du lemme  ci-dessous.
\medskip

{\bf Lemme 3.5.}
{\sl On a des isomorphismes de groupes
 
$$K_0({\bf N}_{A/B})\cong 
K_0({\bf \Lambda})\cong  
K_0({\bf N}_{Cart})\cong
{\cal K}_0({\bf N}_{A/B}).$$}
\medskip

{\sl Dmonstration}. 
On sait que   $K_1({\bf Pic}(A))=A^*$ 
et  $K_0({\bf Pic}(A))=Cl(A)$. 

Par consquent, dans le diagramme \`a lignes exactes
$$\xymatrix{
A^*\ar[r]&
B^*\ar[r]&
K_0({\bf \Lambda})\ar[r]&
Cl(A)\ar[r]&
Cl(B)\cr
A^*\ar[r]\ar[u]^-{\hbox{\sl d\'et}_1}&
B^*\ar[r]\ar[u]^-{\hbox{\sl d\'et}_1}&
K_0({\bf N}_{A/B})\ar[r]\ar[u]^-{\hbox{\bf d\'et}_0}&
Cl(A)\ar[r]\ar[u]^-{\hbox{\sl d\'et}_0}&
Cl(B)\ar[u]^-{\hbox{\sl d\'et}_0}}$$
les  applications $\hbox{\sl d\'et}_1$ et $\hbox{\sl d\'et}_0$ sont  des isomorphismes, ce qui montre que $\hbox{\bf d\'et}_0$ est un isomorphisme, c'est-\`a-dire $K_0({\bf N}_{A/B})\cong K_0({\bf \Lambda})$.

Par ailleurs, puisque 
${\bf i}_A$ et ${\bf i}_B$ sont des \'equivalences de cat\'egories, on a un isomorphime de groupes de $K$-th\'eorie de foncteurs 
$$K_0({\bf \Lambda})\cong K_0({\bf N}_{Cart}).$$

Enfin, d'aprs 2.4, on a un isomorphisme de groupes $K_0({\bf N}_{Cart})\cong
{\cal K}_0({\bf N}_{A/B})$.
Ce lemme 3.5 ach\`eve ainsi la  dmonstration du th\'eor\`eme 3.4. 
\bigskip

{\bf 4. La suite exacte longue de Mayer-Vietoris de $K$-thorie de foncteurs.}
\medskip

{\sl  La suite exacte longue de Mayer Vietoris de $K$-thorie des anneaux.}

\medskip

Nous allons d'abord considrer une suite exacte de Mayer-Vietoris pour la $K$-thorie des anneaux.

\medskip

{\bf D\'efinition 4.1.}  {\sl Un morphisme d'anneaux $i:A\to \widehat A$ est un isomorphisme analytique le long de ${\cal S}$ si

a)  
${\cal S}$ est une partie multiplicative de $A$, compos\'ee d'\'el\'ements centraux et non diviseurs de z\'ero,

b) 
$i({\cal S})$ est une partie multiplicative de $\widehat A$, compos\'ee d'\'el\'ements centraux et non diviseurs de z\'ero,

c)  
pour tout $s$ de ${\cal S}$, l'application $i$ induit un isomorphisme $A/s\cong \widehat A/i(s)$.}
\medskip

On dit alors que $(A,i,\widehat A,{\cal S})$ est un carr\'e de Milnor et on a le diagramme commutatif
$$\xymatrix{
A\ar[r]^{i}\ar[d]&\widehat A\ar[d]\cr
{\cal S}^{-1}A\ar[r]&i({\cal S})^{-1}\widehat A\cr
}$$
\medskip

En $K$-thorie, les carrs de Milnor jouissent de la proprit suivante.

\medskip
{\bf Th\'eor\`eme 4.2.}([K2] p. 400, [W1]) {\sl Soit $(A,i,\widehat A, {\cal S})$ un carr\'e de Milnor. Il existe une suite exacte longue de Mayer-Vietoris en $K$-th\'eorie 
$$\xymatrix{
\cdots \ar[r]^-{\partial}&K_r(A)\ar[r]&K_r({\cal S}^{-1}A)\times K_r(\widehat A)\ar[r]^-{\mu_r}&K_r({\cal S}^{-1}\widehat A)\ar[r]^-{\partial}&\cdots}$$ } 
\medskip

Donnons un exemple :
soit $F$  un corps de nombres d'anneau d'entiers ${\cal O}_F=A$.  L'application $i:A\to\widehat A$ est alors un isomorphisme analytique le long de ${\cal S}={\bf N}\setminus\{0\}$. Dans ce cas, ${\cal S}^{-1}A=F$ et  
$i({\cal S})^{-1}\widehat A$ est l'anneau   ${\bf A}_F$ des adles du corps $F$, restreintes aux places finies. On a par consquent le carr de Milnor $(A,i,\widehat A,{\cal S})$ et le diagramme commutatif
$$\xymatrix{
A\ar[r]^{i}\ar[d]&\widehat A\ar[d]\cr
F\ar[r]&{\bf A}_F\cr
}\leqno (M)$$

\medskip

Appliquons  le foncteur $K_1$ au diagramme $(M)$ ci-dessus. D'aprs la proposition 1.2, on obtient le diagramme
 $$\xymatrix{
A^*\ar[r]\ar[d]&{\cal U}_F\ar[d]\cr
F^*\ar[r]&{\bf J}_F\cr
}\leqno (M_1)$$

En particulier,  pour $r=1$, la suite exacte de Mayer-Vietoris du carr\'e de Milnor $(A,i,\widehat A,{\cal S})$  s'\'ecrit sous la forme 
$$\xymatrix{\cdots K_2({\bf A}_F)\ar[r]^-{\partial}& 
A^*\ar[r]^-{i_1}&
F^*\times {\cal U}_F\ar[r]^-{\mu_1}&
{\bf J}_F\ar[r]^-{\partial}&
K_0(A)\ar[r]^-{i_0}&K_0(F)\times K_0(\widehat A).}$$

L'homomorphisme $i_1$ est  videmment injectif car $A^*$ est contenu dans $F^*$. Par ailleurs, $\ker(i_0)$ s'identifie au groupe des classes $Cl(A)$. On aboutit ainsi \`a la suite exacte

$$\xymatrix{1\ar[r]& 
A^*\ar[r]^-{i}&F^*\times {\cal U}_F\ar[r]^-{\mu}&
{\bf J}_F\ar[r]&Cl(A)\ar[r]& 1}$$
avec $\mu(z,t)=z/t$, ce qui fournit l'isomorphisme classique 
 $${\bf J}_F/F^*\cdot{\cal U}_F \cong Cl(A).$$
 
 Nous verrons en 5.6 une autre application de cette suite exacte longue de Mayer-Vietoris.
\medskip

{\sl La suite exacte longue de Mayer-Vietoris de $K$-thorie de foncteurs.}
\medskip

{\bf Proposition 4.3.}
{\sl Soit ${\cal S}$ la partie multiplicative ${\cal S}={\bf N}\setminus\{0\}$ et soient 
$(A,i_A,\widehat A, {\cal S})$ et $(B,i_B,\widehat B,{\cal S})$ deux carr\'es de Milnor. Soit $\varphi:{\bf Proj}(A)\to{\bf Proj}(B)$ un foncteur additif et
soient
$\widehat\varphi$,
${\cal S}^{-1}\varphi$ et 
${\cal S}^{-1}\widehat\varphi$ les foncteurs induits par $\varphi$. Alors le cube de foncteurs ci-dessous est commutatif \`a isomorphisme pr\`es.

$$\def\objectstyle{\scriptstyle} \def\labelstyle{\scriptstyle}\xymatrix@=20pt@C=20pt{
& {\bf Proj}({\cal S}^{-1}A) \ar[rr]^-{}\ar[dddd]_-{{\cal S}^{-1}\varphi} & & {\bf Proj}({\cal S}^{-1}\widehat A)\ar[dddd]^-{{\cal S}^{-1}\widehat\varphi}
\\
{\bf Proj}(A)\ar[ur]^-{\theta_A}\ar[rr]^-{}\ar[dddd]_-{\varphi} & & {\bf Proj}( \widehat A)\ar[ur]_-{}\ar[dddd]^-{\widehat\varphi}
\\
&
\\
\\
&
{\bf Proj}({\cal S}^{-1}B)\ar[rr]_-{} & & {\bf Proj}( {\cal S}^{-1}\widehat B)
\\
{\bf Proj}(B)\ar[rr]_-{}\ar[ur]^-{}& & {\bf Proj}(\widehat B)\ar[ur]_-{}
}$$
}

{\sl {\sl Dmonstration}}. Les huit foncteurs des faces horizontales sont des foncteurs extensions des scalaires. Par exemple,  $\theta_A$ est l'extension des scalaires associ\'ee au morphisme d'anneaux $A\to {\cal S}^{-1}A$. On a donc
$\theta_A(P)={\cal S}^{-1}A\otimes_AP$ et des formules analogues pour les sept autres foncteurs de ces deux faces. 
\smallskip

Les foncteurs verticaux $\widehat\varphi$, ${\cal S}^{-1}\varphi$ et 
${\cal S}^{-1}\widehat\varphi$ sont d\'efinis \`a partir du foncteur $\varphi$. Puisque $\varphi$ est additif, il existe un $B$-$A$-bimodule $M={}_BM{}_A$ tel que $\varphi=M\otimes_A-$. On pose
\smallskip

\centerline{$\widehat M=\widehat B\otimes_B M\otimes_A \widehat A$ et 
$\widehat\varphi=\widehat M\otimes_{\widehat A}-$,}
\smallskip

\centerline{${\cal S}^{-1}M={\cal S}^{-1}B\otimes_B M\otimes_A {\cal S}^{-1}A$ et ${\cal S}^{-1}\varphi={\cal S}^{-1}M\otimes_{{\cal S}^{-1}A}-$,}
\smallskip

\centerline{${\cal S}^{-1}\widehat M={\cal S}^{-1}\widehat B\otimes_{\widehat B} \widehat M\otimes_{\widehat A}{\cal S}^{-1}\widehat A$ et 
${\cal S}^{-1}\widehat\varphi={\cal S}^{-1}\widehat M\otimes_{{\cal S}^{-1}\widehat A}-$.}
\smallskip

C'est alors un jeu d'\'ecriture fastidieux mais trivial de v\'erifier que le cube est commutatif \`a isomorphisme pr\`es. 

\medskip

{\bf Th\'eor\`eme 4.4.}
{\sl  Soit ${\cal S}$ la partie multiplicative ${\cal S}={\bf N}\setminus\{0\}$ et soient 
$(A,i_A,\widehat A,{\cal S})$ et 
$(B,i_B,\widehat B,{\cal S})$ deux carr\'es de Milnor. Soit $\varphi:{\bf Proj}(A)\to{\bf Proj}(B)$ un foncteur additif. Alors les foncteurs 
$\varphi$, 
$\widehat\varphi_{}$, 
${\cal S}^{-1}\varphi_{}$ et 
${\cal S}^{-1}\widehat\varphi_{}$ 
de la proposition pr\'ec\'edente fournissent une suite exacte longue de Mayer-Vietoris de K-th\'eorie de foncteurs:
$$\xymatrix{
 \cdots \ar[r]& K_{r+1}({\cal S}^{-1}\widehat \varphi_{}) \ar[r]^-{\partial}& K_r(\varphi)\ar[r]& 
 K_r({\cal S}^{-1}\varphi_{})\times K_r(\widehat\varphi_{})\ar[r]^-{\mu_r}& K_r( {\cal S}^{-1}\widehat\varphi_{})\ar[r]^-{\partial}& K_{r-1}(\varphi)\ar[r]&\cdots}$$}

\medskip

{\sl Dmonstration}. Soient ${\cal F}_\varphi$, 
${\cal F}_{{\cal S}^{-1}\varphi}$, ${\cal F}_{\widehat \varphi}$ et ${\cal F}_{{\cal S}^{-1}\widehat \varphi}$ les fibres homotopiques des applications 
$BQ(\varphi)$, $BQ({\cal S}^{-1}\varphi)$, $BQ(\widehat \varphi)$ et $BQ({\cal S}^{-1}\widehat \varphi)$ construites \`a partir du cube ci-dessus. 

Soient ${\cal F}$ et $\widehat {\cal F}$ les fibres homotopiques de ${\cal F}_\varphi\to {\cal F}_{{\cal S}^{-1}\varphi}$ et de
${\cal F}_{\widehat \varphi}\to{\cal F}_{{\cal S}^{-1}\widehat \varphi}$. A partir des applications ${\cal F}_\varphi\to {\cal F}_{\widehat \varphi}$ et
${\cal F}_{{\cal S}^{-1}\varphi}\to {\cal F}_{{\cal S}^{-1}\widehat \varphi}$, on obtient une application $j:{\cal F}\to\widehat {\cal F}$. Les fibrations homotopiques de fibres 
$\cal F$ et $\widehat {\cal F}$ conduisent au diagramme \`a lignes exactes.
$$\xymatrix{
\cdots\ar[r]&
\pi_{r+1}({\cal F})\ar[r]\ar[d]^-{j_{r+1}}&
K_r(\varphi)\ar[r]\ar[d]^-{\iota_r}&
K_r({\cal S}^{-1}\varphi)\ar[r]\ar[d]^-{\overline\iota_r}&
\pi_r({\cal F})\ar[r]\ar[d]^-{j_r}&
K_{r-1}(\varphi)\ar[r]\ar[d]^-{\iota_{r-1}}&
\cdots\cr
\cdots\ar[r]&
\pi_{r+1}(\widehat {\cal F})\ar[r]&
K_r(\widehat \varphi)\ar[r]&
K_r({\cal S}^{-1}\widehat \varphi)\ar[r]&
\pi_{r}( \widehat {\cal F})\ar[r]&
K_{r-1}(\widehat \varphi)\ar[r]&
\cdots
}$$
Pour obtenir la suite exacte de Mayer-Vietoris propos\'ee, montrons que $j_r$ est un isomorphisme pour tout $r\geq 0$.  D'apr\`es Quillen  et Gersten (([Q], [G]), la fibre homotopique de $BQ{\bf Proj}(A)\to BQ{\bf Proj}({\cal S}^{-1}A)$  a le type d'homotopie de l'espace classi\-fiant $B{\bf T}^1_{\cal S}(A)$ o\`u ${\bf T}^1_{\cal S}(A)$ est la cat\'egorie des $A$-modules de ${\cal S}$-torsion ayant une r\'esolution de longueur inf\'erieure \`a $1$ par des $A$-modules projectifs de type fini.

La fibre homotopique de $B{\bf T}^1_{\cal S}(A)\to B{\bf T}^1_{\cal S}(B)$ a le type d'homotopie de ${\cal F}$ et celle de 
$B{\bf T}^1_{\cal S}(\widehat A)\to B{\bf T}^1_{\cal S}(\widehat B)$ a le type d'homotopie de $\widehat {\cal F}$. De plus les applications 
$B{\bf T}^1_{\cal S}(A)\to B{\bf T}^1_{\cal S}(\widehat A)$ et $B{\bf T}^1_{\cal S}(B)\to B{\bf T}^1_{\cal S}(\widehat B)$  induisent une application 
${\cal F}\to\widehat {\cal F}$  qui est homotope \`a l'application $j$. 

Les suites exactes longues des fibrations homotopiques $B{\bf T}^1_{\cal S}(A)\to B{\bf T}^1_{\cal S}(B)$ et $B{\bf T}^1_{\cal S}(\widehat A)\to B{\bf T}^1_{\cal S}(\widehat B)$ conduisent au diagramme commutatif \`a lignes exactes
$$\xymatrix{
\cdots\ar[r]&
\pi_{r+1}({\cal F})\ar[r]\ar[d]^-{j_{r+1}}&
\pi_{r+1}(B{\bf T}^1_{\cal S}(A))\ar[r]\ar[d]^-{i_{r+1}}&
\pi_{r+1}(B{\bf T}^1_{\cal S}(B))\ar[r]\ar[d]^-{i'_{r+1}}&
\pi_r({\cal F})\ar[r]\ar[d]^-{j_{r}}&\cdots\cr
\cdots\ar[r]&
\pi_{r+1}(\widehat {\cal F})\ar[r]&
\pi_{r+1}(B{\bf T}^1_{\cal S}(\widehat A))\ar[r]&
\pi_{r+1}(B{\bf T}^1_{\cal S}(\widehat B))\ar[r]&
\pi_r(\widehat {\cal F})\ar[r]&\cdots
}$$

Or $i:A\to \widehat A$ et $j:B\to \widehat B$ sont des isomorphismes analytiques le long de ${\cal S}$. D'apr\`es [K2], p. 400-402, on sait que dans ce cas, les foncteurs 
 $\widehat A\otimes_A -:{\bf T}^1_{\cal S}(A)\to {\bf T}^1_{\cal S}(\widehat A)$ et $\widehat B\otimes_B -:{\bf T}^1_{\cal S}(B)\to {\bf T}^1_{\cal S}(\widehat B)$ sont des \'equivalences de cat\'egories. Il en r\'esulte que $i_r$ et $i'_r$ sont des isomorphismes pour tout $r\geq 0$. Par le lemme des cinq, on en d\'eduit que $j_r$ est un isomorphisme pour tout $r\geq 0$.
 
Ceci permet de d\'efinir le connectant de la suite exacte de Mayer-Vietoris de $K$-th\'eorie de foncteurs comme la composition des homomorphismes 
$$\xymatrix{
K_r({\cal S}^{-1}\widehat \varphi)\ar[r]&
K_r(\widehat {\cal F})\ar[r]^-{j_r^{-1}}&
\pi_r({\cal F})\ar[r]&
K_{r-1}(\varphi)}.$$

 \bigskip

{\bf  5. Cas des corps de nombres et du foncteur norme.}
\medskip

Soit $F/L$ une extension  de corps de nombres, d'anneaux d'entiers $A={\cal O}_F$, $B={\cal O}_L$.  
Le morphisme d'anneaux   $i_A:A\to \widehat A$ est un isomorphisme analytique le long de $S={\bf N}\setminus\{0\}$. On a  ${\cal S}^{-1}A=F$ et 
${\cal S}^{-1}\widehat A={\bf A}_F$.
\medskip

Soit $\varphi:{\bf Proj}(A)\to{\bf Proj}(B)$ un foncteur additif. Pour harmoniser les notations, notons 
$\varphi_{\widehat A/\widehat B}=\widehat \varphi$,  
\ $\varphi_{F/L}={\cal S}^{-1}\varphi$ et  
$\varphi_{{\bf A}_F/{\bf A}_L}={\cal S}^{-1}\widehat \varphi$
les foncteurs introduits en 4.3, c'est-\`a-dire :
 
$\varphi_{\widehat A/\widehat B}:{\bf Proj}(\widehat A)\to {\bf Proj}(\widehat B)$,

$\varphi_{F/L}:{\bf Proj}(F)\to {\bf Proj}(L)$

 et
$\varphi_{{\bf A}_F/{\bf A}_L}:{\bf Proj}({\bf A}_F)\to 
{\bf Proj}({\bf A}_L)$.

\medskip

La suite exacte longue de Mayer-Vietoris de $K$-th\'eorie de foncteurs 4.4 s'\'ecrit 
$$\xymatrix{
\cdots \ar[r]& 
K_{1}(\varphi_{{\bf A}_F/{\bf A}_L)} \ar[r]^-{\partial}& 
K_0(\varphi)\ar[r]& 
K_0(\varphi_{F/L})\times 
K_0(\varphi_{\widehat A/\widehat B})\ar[r]^-{\mu_0}&
K_0(\varphi_{{\bf A}_F/{\bf A}_L}).
}$$

\medskip

Convenons galement des allgements  de notations suivants.

${}_\varphi K_r(A)=\ker K_r(A)\to K_r(B)$,\ \  
$K_r(A)\varphi=\hbox{coker\ } K_r(A)\to K_r(B)$,

${}_\varphi K_r(F)=\ker K_r(F)\to K_r(L)$,\ \  
$K_r(F)\varphi=\hbox{coker\ } K_r(F)\to K_r(L)$,

${}_\varphi K_r(\widehat A)=\ker K_r(\widehat A)\to K_r(\widehat B)$,\ \  
$K_r(\widehat A)\varphi=\hbox{coker\ } K_r(\widehat A)\to K_r(\widehat B)$,

${}_\varphi K_r({\bf A}_F)=\ker K_r({\bf A}_F)\to K_r({\bf A}_L)$,\ \  
$K_r(A)\varphi=\hbox{coker\ } K_r({\bf A}_F)\to K_r({\bf A}_L)$,
\medskip

Cette  suite exacte de Mayer-Vietoris de $K$-th\'eorie de foncteurs peut \^etre  pr\'ecis\'ee, notamment en degrs $0$ et $1$. Pour cela rappelons le rsultat suivant, d\^u  C. Soul.

\medskip

{\bf Thorme 5.1.} ([Sou])
{\sl Soit $F$ un corps de nombres d'anneaux d'entiers $A$. Alors pour $r>0$, l'application 
$K_r(A)\to K_r(F)$ est injective.}
\medskip

Ce rsultat nous permet de dcouper la suite de Mayer-Vietoris en morceaux sous la forme suivante.
\medskip

{\bf Proposition 5.2.}
{\sl Soit $F/L$ une extension de corps de nombres, d'anneaux d'entiers 
$A={\cal O}_F$ et $B={\cal O}_L$ . 
Soit $\varphi:{\bf Proj}(A)\to{\bf Proj}(B)$ un foncteur additif. Pour $r\geq 0$, on a la suite exacte

$$\xymatrix{
&1\ar[r]&
{}_\varphi K_{r+1}(A)\ar[r]^-{\Delta}&
{}_\varphi K_{r+1}(F)\times{}_\varphi K_{r+1}(\widehat A)\ar[r]^-{\mu_{r+1}}&
{}_\varphi K_{r+1}({\bf A}_F)\ar[dldl]^-{\partial}&\cr
\cr
&&K_r(\varphi)\ar[r]^-{i_r}&
K_r(\varphi_{F/L})\times K_r(\varphi_{\widehat A/\widehat B})\ar[r]^-{\mu_r}&
K_r(\varphi_{{\bf A}_F/{\bf A}_L})}$$
} 

\medskip
  
{\sl Dmonstration}. Considrons  le diagramme commutatif suivant dans lequel les lignes et les colonnes sont extraites des suites exactes   de Mayer-Vietoris.

 $$\def\objectstyle{\scriptstyle} \def\labelstyle{\scriptstyle}\xymatrix@=10pt@C=10pt{
&
\vdots\ar[d]&
{}_\varphi K_{r+1}(A)\ar[dd]^(.3){\Delta}&
\vdots\ar[d]&
\vdots\ar[d]
\cr
\cdots\ar[r]\ar[dd]&
K_{r+1}(\varphi_{F/L})\times 
K_{r+1}(\varphi_{\widehat A/\widehat B)}\ar[dd]\ar[rr]&&
K_{r+1}(F)\times K_{r+1}(\widehat A)\ar[r]^-{(\varphi_{F}\times \varphi_{\widehat A})_{{}_{r+1}}}\ar[dd]^-{\widehat\mu_{r+1}}&
K_{r+1}(L)\times K_{r+1}(\widehat B)\ar[dd]
\cr
&&
{}_\varphi K_{r+1}(F)\times{}_\varphi K_{r+1}(\widehat A)\ar[ru]\ar[dd]^(.7){\mu_{r+1}}
&&&
\cr
\cdots 
K_{r+2}({\bf A}_L)\ar[r]\ar[dd]^-{\partial_B}&
K_{r+1}(\varphi_{{\bf A}_F/{\bf A}_L})\ar[rr]^(.3){\widehat \iota}\ar[dd]^{\partial'}\ar[dr]^(.6){\iota_{r+1}}&&
K_{r+1}({\bf A}_F)\ar[r]^-{(\varphi_{{\bf A}_F})_{{}_{r+1}}}\ar[dd]&
K_{r+1}({\bf A}_L)\ar[dd]
\cr
&&
{}_\varphi K_{r+1}({\bf A}_F)\ar[ru]\ar[dl]_(.3){\partial}&&
\cr
\cdots
K_{r+1}(B)\ar[r]\ar[dd]^-{j_{r+1}}&
K_r(\varphi)\ar[dd]^-{i_r}\ar[rr]&&
K_r(A)\ar[r]^-{\varphi_{{}_r}}&
K_r(B)
\cr
&\cr
\cdots K_{r+1}(L)\times K_{r+1}(\widehat B)\ar[r]&
K_r(\varphi_{F/L})\times K_r( \varphi_{\widehat A/\widehat B})\ar[dd]^-{\mu_r}&&&
\cr
&\cr
&K_r(\varphi_{{\bf A}_F})&&&
}$$

Pour $r\geq 0$, l'homomorphisme $K_{r+1}(B)\to K_{r+1}(L)$ est injectif (Soul).  Par consquent, l'homomorphisme   $j_{r+1}:K_r(B)\to K_r(L)\times K_r(\widehat B)$ est galement injectif et  $\partial_B=0$. Par chasse dans le diagramme, on construit un homomorphisme 
$\partial:{}_\varphi K_{r+1}({\bf A}_F)\to K_r(\varphi)$ 
de m\^eme image que $\partial'$. Soit $\mu_{r+1}$ la restriction de $\widehat\mu_{r+1}$ \`a 
${}_\varphi K_{r+1}(F)\times\ {}_\varphi K_{r+1}(\widehat A)$. 

On a 
$$\ker\partial=\iota_{r+1}(\ker\partial')=
\mu_{r+1}({}_\varphi K_{r+1}(F)\times{}_\varphi K_{r+1}(\widehat A))=
{}_\varphi K_{r+1}(F)\cdot{}_\varphi K_{r+1}(\widehat A)=\hbox{im\ }(\mu_{r+1}).$$

Notons 
$\Delta:{}_\varphi K_{r+1}(A)\to
{}_\varphi K_{r+1}(F)\times{}_\varphi K_{r+1}(\widehat A)$ l'homomorphisme d\'efini par 
$\Delta (z)=(z,z)$, en identifiant 
${}_\varphi K_{r+1}(A)$  un sous-groupe de 
${}_\varphi K_{r+1}(F)$ et de 
${}_\varphi K_{r+1}(\widehat A)$. Il est clair que $\ker\mu_{r+1}=\hbox{im\ }(\Delta).$
Ceci ach\`eve la dmonstration.

\bigskip

{\sl La suite exacte de Mayer-Vietoris du foncteur norme.}

\medskip

Dtaillons la suite exacte 5.2 pour le foncteur norme  ${\bf N}_{A/B}$ et pour $r=0$.
\medskip

{\bf Proposition 5.3.}
{\sl 
Soit  $F/L$ une extension  de corps de nombres, d'anneaux d'entiers $A={\cal O}_F$ et $B={\cal O}_L$. 
Les   groupes  de $K$-th\'eorie des foncteurs 

${\bf N}_{F/L}:{\bf Proj}(F)\to {\bf Proj}(L),$

${\bf N}_{\widehat A/\widehat B}:{\bf Proj}(\widehat A)\to {\bf Proj}(\widehat B)$ et

${\bf N}_{{\bf A}_F/{\bf A}_L}:{\bf Proj}({\bf A}_F)\to {\bf Proj}({\bf A}_L)$
induits par le foncteur norme ${\bf N}_{A/B}$   
 sont dtermins ainsi 

$$K_0({\bf N}_{F/L})=L^*_N,$$

$$K_0({\bf N}_{\widehat A/\widehat B})=
({\cal U}_L)_N,$$

$$K_0({\bf N}_{{\bf  A}_F/{\bf A}_L})=
\left({\bf J}_{L}\right)_N.$$
}
\medskip

{\sl Dmonstration.}  
Notons toujours $N:F^*\to L^*$ la norme. La suite exacte de Bass 2.2 du foncteur ${\bf N}_{F/L}$  montre immdiatement qu'on a 
$K_0({\bf N}_{F/L})=L^*_N$.
\smallskip

La suite exacte de Bass du foncteur 
${\bf N}_{\widehat A/\widehat B}$ s'crit:

$$\xymatrix{
K_1(\widehat A)\ar[r]^{N}& 
K_1(\widehat B)\ar[r]&
K_0({\bf N}_{{\widehat A/\widehat B}})\ar[r]&
K_0(\widehat A)\ar[r]^N&
 K_0(\widehat B).}$$
D'aprs 1.1,  $K_1(\widehat A)={\cal U}_F$ et 
$K_1(\widehat B)={\cal U}_L$. Pour montrer la relation 
 $K_0({\bf N}_{\widehat A/\widehat B})=({\cal U}_L)_N$, il suffit de montrer que la norme 
$N:K_0(\widehat A)\to K_0(\widehat B)$ est injective. 
Mais d'aprs 1.3, 
$K_0(\widehat A)$ s'identifie au produit born $\prod^{(b)}_{w\in V(F)}K_0(A_w)$. Ce dernier est un sous groupe de 
$\prod_{w\in V(F)}K_0(A_w)={\bf Z}^{V(F)}$.   De mme, 
$K_0(\widehat B)$ est un sous-groupe de 
$\prod_{v\in V(L)}K_0(B_v)={\bf Z}^{V(L)}$. Dans le diagramme commutatif 
$$\xymatrix{
K_0(\widehat A)\ar[r]^N\ar[d]&
K_0(\widehat B)\ar[d]\cr
{\bf Z}^{V(F)}\ar[r]^{N'}&
{\bf Z}^{V(L)} 
}$$
les flches verticales et la flche $N'$ sont des injections. Il en rsulte que la norme $N=N_{\widehat A/\widehat B}$ est galement injective. On en dduit l'galit 
$K_0({\bf N}_{\widehat A/\widehat B})=
({\cal U}_L)_N.$ 
 \smallskip
 
Pour montrer la dernire galit 
$K_0({\bf N}_{{\bf  A}_F/{\bf A}_L})=
\left({\bf J}_{L}\right)_N$, on considre la suite exacte 
$$\xymatrix{
K_1({\bf A}_F)\ar[r]^N& 
K_1({\bf A}_L)\ar[r]&
K_0(N_{{\bf A}_F/{\bf  A}_L)}\ar[r]&
K_0({\bf A}_F)\ar[r]^N&
 K_0({\bf A}_L).}$$

Le calcul 1.3 et un raisonnement analogue  celui ci-dessus montre que  la norme 
$N:K_0({\bf A}_F)\to K_0({\bf A}_L)$ est galement injective. La formule $K_0({\bf N}_{{\bf  A}_F/{\bf A}_L})=
\left({\bf J}_{L}\right)_N$ en rsulte.

\medskip

Pour d\'ecrire la suite exacte de Mayer-Vietoris 5.2 du foncteur
 ${\bf N}_{A/B}$, gardons les conventions dj employes, o on note indiffremment $N:X\to Y$ la norme entre les groupes correspondants. On a donc 
${}_NX=\ker N$ et 
$Y_N=\hbox{coker} N.$
\bigskip

{\bf Thorme  5.4} (suite exacte des noyaux-conoyaux)

{\sl Soit $F/L$ une extension  de corps de nombres. On a une suite exacte  
$$\xymatrix{
&1\ar[r]&
{}_NA^*\ar[r]^-{\Delta}&
{}_NF^*\times{}_N{\cal U}_F\ar[r]^-{\mu_1}&
{}_N{\bf J}_F\ar[dldl]^-{\partial}&\cr
\cr
&&K_0({\bf N}_{A/B})\ar[r]^-{i}&
L^*_N\times  ({\cal U}_L)_N\ar[r]^-{\mu}&
({\bf J}_L)_N\ar[r]^-{\partial'}&
Cl(B)_N\ar[r]&1.}$$} 

\medskip

{\sl Dmonstration.}
Le rsultat 5.2 fournit la suite exacte jusqu'au terme $(J_L)_N$. Pour l'exactitude en $({J_L})_N$ et la surjectivit de $\partial'$, on procde au dla\c cage de la $K$-thorie algbrique ([B], [K1], [Wa]).  La suite exacte longue de Mayer-Vietoris de $K$-thorie des foncteurs se prolonge aux entiers ngatifs et dbute par

$$\xymatrix{
\ar[r]&
K_0({\bf N}_{{\bf A}_F/{\bf A}_L})=(J_L)_N\ar[r]^-{\partial'}&
K_{-1}({\bf N}_{A/B})=Cl(B)_N\ar[r]&
K_{-1}({\bf N}_{F/L})\times 
K_{-1}({\bf N}_{\widehat A/\widehat B})=0\ar[r]&\cdots
}$$
ce qui montre la surjectivit de $\partial'$.
\bigskip
\vfill\eject

{\bf Corollaire 5.5.} {\sl Soit $F/L$ une extension  de corps de nombres. 

a) Le groupe $K_0({\bf N}_{A/B})$ s'ins\`ere dans les deux suites exactes ci-dessous.

$$\xymatrix{
1\ar[r]&
B^*_N\ar[r]^-{\sigma}&
K_0({\bf N}_{A/B})\ar[r]^-{\rho}\ar@{=}[d]&
{}_NCl(A)\ar[r]&
1\cr
&
{}_N{\bf J}_F\ar[r]^-{\partial}&
K_0({\bf N}_{A/B})\ar[r]^-{i}&
L^*_N\times 
({\cal U}_L)_N\ar[rr]^-{\mu}&&
({\bf J}_L)_N\cr
}$$

b) Posons 
$\alpha=\rho\partial$ et $\beta=i\sigma$. On a un  isomorphisme de groupes 
 $$coker\ \alpha\cong\ker\mu/\hbox{im\ }\beta.\leqno (5.5)$$}

\medskip

Dcrivons  les morphismes $\partial$, $i$ et $\mu$ du diagramme 5.5 (les morphismes $\sigma$ et $\rho$ ont t d\'ecrits en 3.4).

 \medskip
 Pour d\'ecrire $\partial$ et $i$, on identifie  $K_0({\bf N}_{A/B})$ avec le groupe ${\cal K}_0({\bf N}_{A/B})$ du th\'eor\`eme 3.4.  Soit $z=(z_w)_{w\in V(F)}$ un \'el\'ement de ${\bf J}_F$, avec $z_w$ dans $F^*_w$ de valuation $r_w\in {\bf Z}$.
 Soit $\mathfrak p_w$ l'id\'eal de la valuation $w$. Introduisons l'id\'eal fractionnaire $I_z$ de $F$ d\'efini par 
 $I_z=\prod_{w\in V(F)}\mathfrak p_{w}^{r_w}$. La relation $N(I_z)=I_{N(z)}$ montre que si 
 $z$ appartient \`a ${}_N{\bf J}_F$, alors $N(I_z)$ est un \'el\'ement de $B^*$. On en d\'eduit  que $[1_B, I_z]$ appartient
 \`a ${\cal K}_0({\bf N}_{A/B})$ et que $$\partial(z)=[1_B,I_z].$$
 \medskip
 
 Soit $[t,I]$ un \'el\'ement de ${\cal K}_0({\bf N}_{A/B})$. Posons $[t]=t\hbox{\ mod\ }N(F^*)$, donc $[t]\in  L^*_N$. Soit $v$ une place finie de $L$, d'uniformisante $\pi_{v}\in L_v$. Notons $r\in{\bf Z}$  la valuation de 
 $t\in L_v$ et introduisons l'unit\'e $u_v(t)=t\pi^{-r}$ de $B^*_v$. On pose 
 $[u_v(t)]=u_v(t)\hbox{\ mod\ }N(A^*_w)$, o\`u $w$ est une place de $A$ au-dessus de $v$. Enfin  
 $[u(t)]$ d\'esigne la famille $([u_v(t)])_{v\in V(L)}$.
 On a alors $$i([t,I])=([t],[u(t)]).$$
 
 \medskip
  Pour $([t],[u])$ dans $L^*_N\times ({\cal U}_L)_N$, on a enfin
  $$\mu([t],[u])=([t/u_v])_{v\in V(L)}.$$
  
 \medskip

 {\bf Remarque 5.6.} (J. Berrick, M. Lim)
Dans le diagramme ci-dessous, les lignes sont exactes car extraites des suites de Mayer-Vietoris de $K$-thorie d'anneaux.

$$\xymatrix{
1\ar[r]&B^*\ar[r]&L^*\times {\cal U}_L\ar[r]&{\bf J}_L\ar[r]&Cl(B)\ar[r]&1\cr
1\ar[r]&A^*\ar[r]\ar[u]^{\phi_1}&F^*\times {\cal U}_F\ar[r]\ar[u]^{\phi_1}&
{\bf J}_F\ar[r]\ar[u]^{\phi_1}&Cl(A)\ar[r]\ar[u]^{\phi_0}&1\cr
}$$

 On note indiffremment $\phi_1$   les homomorphismes induits en $K$-thorie par les foncteurs $\varphi$, $\varphi_{F/L}$, $\varphi_{\widehat A/\widehat B}$ et
 $\varphi_{{\bf A}_F/{\bf A}_L}$ sur le groupe  $K_1$.
 En considrant ce diagramme comme un complexe double de premier quadrant, la suite spectrale cohomologique  associe  ce complexe converge vers $0$.  Pour $0\leq p\leq 3$, on en dduit des isomorphismes 
 $$E^2_{p,1}\cong E^2_{p+2,0}.$$
 \medskip

Dans le cas du foncteur norme, c'est--dire si 
$\varphi={\bf N}_{A/B}$, l'isomorphisme $E^2_{2,1}\cong E^2_{4,0}$  est l'isomorphisme $coker(\alpha)\cong\ker(\mu)/\hbox{im}(\beta)$ de 5.5.
 \medskip

 Pour $r>1$, compte tenu du thorme de Soul [Sou], un analogue   du diagramme  ci-dessus est le suivant :
 $$\xymatrix{
1\ar[r]&
K_r(B)\ar[r]&
K_r(L)\times K_r(\widehat B)\ar[r]&
K_r({\bf A}_L)\ar[r]&
1\cr
1\ar[r]&
K_r(A)\ar[r]\ar[u]^{\phi_r}&
K_r(F)\times K_r(\widehat A)\ar[r]\ar[u]^{\phi_r}&
K_r({\bf A}_F)\ar[r]\ar[u]^{\phi_r}
&1\cr
}$$
Le m\^eme raisonnement conduit  des isomorphimes 
 $E^2_{1,2}\cong E^2_{3,0}$ et $E^2_{2,1}\cong 0$.
 \bigskip

{\sl Cas des extensions galoisiennes.}
\medskip

Dans le cas d'une extension galoisienne de groupe de Galois $G$,  on a des isomorphismes

$$K_0({\bf N}_{F/L})\cong \widehat H^0(G,F^*),$$

$$K_0({\bf N}_{\widehat A/\widehat B})\cong 
\widehat H^0(G,{\cal U}_F)$$
et
$$K_0({\bf N}_{{\bf A}_F/{\bf A}_L})\cong 
\widehat H^0(G,{\bf J}_F).$$

La suite exacte 5.4 s'crit donc ainsi :
$$\xymatrix{
&1\ar[r]&
{}_NA^*\ar[r]^-{\Delta}&
{}_NF^*\times{}_N{\cal U}_F\ar[r]^-{\mu_1}&
{}_N{\bf J}_F\ar[dldl]^-{\partial}&\cr
\cr
&&K_0({\bf N}_{A/B})\ar[r]^-{i}&
\widehat H^0(G, F^*)\times\widehat H^0(G,{\cal U}_F)  \ar[r]^-{\mu}&
\widehat H^0(G,{\bf J}_F)\ar[r]^-{\partial'}&
Cl(B)_N\ar[r]&1.}$$  
 
\medskip

  L'algbre homologique va nous permettre de pousser plus avant les calculs. Fixons pour cela quelques notations.

Notons
$$h_0:\widehat H^0(G, F^*)\to 
\widehat H^0(G, {\bf J}_F)$$ 
et
$$j_0:\widehat H^0(G, {\cal U}_F)\to 
\widehat H^0(G, {\bf J}_F)$$

les morphismes de groupes  induits par les inclusions canoniques 
$F^*\to {\bf J}_F$ et ${\cal U}_F\to {\bf J}_F$.

Pour $(x,y)\in 
\widehat H^0(G, F^*)\times \widehat H^0(G,{\cal U}_F)$, on a donc $\mu(x,y)=h_0(x)j_0(y)^{-1}$.

\medskip

Dans le paragraphe suivant, nous aurons besoin de la description prcise du sous-groupe $\hbox{im\ }j_0\ \cap \ \hbox{im\ }h_0$ de 
$\widehat H^0(G, {\bf J}_F)$. Nous avons rassembl dans la proposition 5.7 divers rsultats sans doute bien connus sur ce groupe. Pour la commodit du lecteur, nous en proposons des dmonstrations.

 \medskip
 {\bf Proposition 5.7.}
 {\sl Soit $F/L$ une extension galoisienne de corps de nombres, de groupe de Galois $G$.
 \smallskip

 1) L'homomorphisme $j_0$ ci-dessus est injectif.
 
 \smallskip
 
 2) On a une suite exacte 
 $$\xymatrix{
 1\ar[r]&
 \hbox{im\ }h_0\ \cap\  \hbox{im\ }j_0\ar[r]&
 \hbox{im\ }j_0\ar[r]^-{}&
 G_{ab}.
 }$$
 
 \smallskip
 
 3) Lorsque l'extension 
est cyclique, l'homomorphisme $h_0$ est injectif.
 
 \smallskip
 
 4) Lorsque l'extension  est cyclique ramifie de degr premier $\ell$, 
 $\hbox {im \ }h_0\cap\hbox{im\ }j_0$ est un 
 ${\bf F}_\ell$-espace vectoriel de dimension $t-1$, o $t$ est le nombre de places finies de $L$ ramifies dans  $F$.
}
\medskip

{\sl Dmonstration.}

1) Injectivit de $j_0$.
Soient $v$ une place de $L$ et $w$ une place de $F$ au-dessus de $v$. On pose $G_w=Gal(F_w/L_v)$. On dsigne par $V(L)$ l'ensemble des places finies de $L$ et par
 ${\cal R}(F/L)$ l'ensemble des places finies de $L$ ramifies dans $F$.

 Montrons en premier lieu qu'on a des isomorphismes de groupes 
 $$ \widehat H^0(G, {\cal U}_F)\cong 
 \oplus_{v\in{\cal R}(F/L)}\widehat H^0(G_w,A^*_w)$$
 et
 $$\widehat H^0(G, {\bf J}_F)\cong
 \oplus_{v\in{V(L)}}\widehat H^0(G_w,F^*_w).$$
 
 \smallskip

 On d\'ecompose le module galoisien 
${\cal U}_F$ en produit direct de modules galoisiens
$${\cal U}_F=\prod_{v\in V(L)}A^*_v$$ avec 
$$A^*_v=\prod_{w\in V(F),w\mid v}A^*_w$$ d'o\`u
$$\widehat H^0(G,{\cal U}_F)=\prod_{v\in V(L)}\widehat H^0(G,A^*_v).$$
Mais $A^*_v=\hbox{Ind\ }_G^{G_v}$ est un module induit et d'apr\`es le lemme de Shapiro, on a 
$$\widehat H^0(G,A^*_v)
=\widehat H^0(G_v,A^*_w).$$ 

Pour conclure, on remarque que $\widehat H^0(G_v,A^*_w)=B^*_v/N(A^*_w)$ et que si $v$ est non ramifi\'ee dans $F$, $B^*_v=N(A^*_w)$ ([N], V.1.2, p. 319). Pour $v$ non ramifie, on a ainsi 
$\widehat H^0(G, A^*_v)=0$. Il reste donc 
$$\widehat H^0(G,{\cal U}_F)=\oplus_{v\in{\cal R}(F/L)}\widehat H^0(G_v,A^*_w).$$ 
\medskip

Le calcul de $\widehat H^0(G,{\bf J}_F)$ s'effectue de manire analogue, par passage  la limite inductive sur l'ensemble des places finies de $L$.
Introduisons  les modules galoisiens 
$F^*_v=\prod_{w\mid v}F^*_w$ et $A^*_v=\prod_{w\mid v} A_w^*.$ Soit $S$ un ensemble fini de places de $L$. On pose  
${\bf J}_F^S=\prod_{v\in S} F^*_v\times\prod_{v\not\in S}A^*_v$ de sorte que ${\bf J}_F=\displaystyle\lim_{\to S} {\bf J}^S_F$.

 \smallskip
 
 On a par consquent
 $\widehat H^0(G, {\bf J}_F)=
 \displaystyle\lim_{\to S}\widehat H^0(G, {\bf J}^S_F)$.
 
 Or 
 $$\widehat H^0(G, {\bf J}^S_F)=
 \oplus_{v\in S}\widehat H^0(G, F_v^*)\times
 \oplus_{v\not\in S}\widehat H^0(G, A^*_v).$$
 
 Par ailleurs, si $v$ est une place  de $L$ non ramifie dans $F$, on a vu que 
 $\widehat H^0(G, A^*_v)=0$. L'ensemble des places ramifies tant fini, on a ${\cal R}(F/L)\subset  S$ pour $S$ assez grand et dans ce cas, on a 
 $$\widehat H^0(G, {\bf J}^S_F)=\oplus _{v\in S}\widehat H^0(G, F^*_v),$$
 d'o 
 $$\widehat H^0(G, {\bf J}_F)=\oplus_{v\in V(L)}\widehat H^0(G, F_v^*).$$
 
 L'emploi du lemme de Shapiro conduit   l'isomorphisme
 $\widehat H^0(G, F^*_v)=\widehat H^0(G_w, F^*_w)$, ce qui donne finalement
 $$\widehat H^0(G, {\bf J}_F)\cong
 \oplus_{v\in{V(L)}}\widehat H^0(G_w,F^*_w).$$

\smallskip

Montrons  prsent l'injectivit de l'homomorphisme $j_0$. D'aprs les calculs qui prcdent, il suffit  de vrifier que pour toute place $v$ de $L$ et toute place $w$ de $F$ au-dessus de $v$, l'homomorphisme 
$$\widehat H^0(G_w,A^*_w)\to \widehat H^0(G_w,F^* _w)$$
est injectif. La suite exacte courte de valuation
$$1\to A^*_w\to F^*_w\to{\bf Z}\to 0$$
conduit  la suite exacte
$$\widehat H^{-1}(G_w,{\bf Z})\to
\widehat H^0(G_w,A^*_w)\to \widehat H^0(G_w,F^*_w).$$

L'opration de $G_w$ sur ${\bf Z}$ tant triviale, le groupe 
$\widehat H^{-1}(G_w,{\bf Z})$ est un quotient du noyau de la multiplication par $\# G_w$ dans ${\bf Z}$. Ce noyau est videmment nul donc 
$\widehat H^{-1}(G_w,{\bf Z})=0$, ce qui fournit l'injectivit recherche.

\smallskip

2) Dmonstration de la suite exacte 
$$\xymatrix{
 1\ar[r]&
 \hbox{im\ }h_0\ \cap\  \hbox{im\ }j_0\ar[r]&
 \hbox{im\ }j_0\ar[r]^-{}&
 G_{ab}.
 }$$
 
 Rappelons que dans ce texte, nous dsignons par ${\bf J}_F$ le groupe des idles aux places finies du corps $F$. Notons ${\bf J}_F'$ le groupe des idles de $F$  toutes les places de $F$ (y compris les places  l'infini). Soit $r$ le nombre de places infinies de $L$ ramifies dans $F$. On a 
 $$\widehat H^0(G, {\bf J}_F')=
 {\bf C}_2^r\times \widehat H^0(G, {\bf J}_F).$$
 
 Introduisons les applications 
 
 $$h'_0:\widehat H^0(G, F^*)\to \widehat H^0(G, {\bf J}_F')$$  
 
 $$j'_0:\widehat H^0(G,{\cal U}_F)\to \widehat H^0(G, {\bf J}_F')$$
 
 et
 
$$i:\widehat H^0(G, {\bf J}_F)\to \widehat H^0(G, {\bf J}_F').$$

On a 
$i\circ h_0=h'_0$ et $i\circ j_0=j'_0$.

La suite exacte courte 
$$1\to F^*\to {\bf J}_F'\to C_F\to 1$$ fournit la suite exacte  $$\xymatrix{
\widehat H^{-1}(G, C_F)\ar[r]&
\widehat H^0(G,F^*)\ar[r]^{h'_0}&
\widehat H^0(G,{\bf J}_F')\ar[r]^{s_0}&
\widehat H^0(G, C_F)\ar[r]&
\widehat H^1(G, F^*).
}$$
D'aprs le thorme de Tate-Nakayama et la thorie globale du corps de classes, le cup-produit fournit un isomorphisme 
$\widehat H^0(G, C_F)\cong \widehat H^{-2}(G,{\bf Z})$ et ce dernier groupe est gal  
$H_1(G,{\bf Z})=G_{ab}=G/[G,G]$. Par ailleurs,  d'aprs le thorme 90 de Hilbert, $H^1(G, F^*)=0$. On a donc la suite exacte
$$\xymatrix{
\widehat H^0(G, F^*)\ar[r]^{h'_0}&
\widehat H^0(G, {\bf J}_F')\ar[r]^{s_0}&
G_{ab}\ar[r]&
1.}$$ 
L'application compose 
$$\varphi_0:\widehat H^0(G,{\cal U}_F)\to \widehat H^0(G,C_F)$$ dfinie par 
$\varphi_0=s_0\circ j'_0$ est donc  valeurs dans $G_{ab}$ et de noyau isomorphe  $$\ker s_0\cap\hbox{im\ }j'_0=
\hbox{im\ }h'_0\cap\hbox{im\ }j'_0.$$
Mais $\hbox{im\ }h'_0\cong \hbox{im\ }h_0$ et 
$\hbox{im\ }j'_0\cong \hbox{im\ }j_0$
car $i$ est injective, d'o le rsultat.

\smallskip

3) Injectivit de $h_0$. C'est le principe de Hasse ([N-S-W], p. 375)

\smallskip

4) Dtermination de $\hbox{im\ }h_0\ \cap\  \hbox{im\ }j_0.$

Posons $\varphi_r=s_r\circ j_r$ avec 
$$j_r:\widehat H^r(G,{\cal U}_F)\to \widehat H^r(G,{\bf J}'_F)$$
et 
$$s_r:\widehat H^r(G,{\bf J}'_F)\to \widehat H^r(G,C_F)$$
et donnons l'expression de $\varphi_2$. D'aprs le lemme de Shapiro, on a les galits 
$$\widehat H^2(G,{\cal U}_F)\cong
\oplus_{v\in {\cal R}(F/L)}H^2(G_w,A^*_w)$$ et 
$$\widehat H^2(G,{\bf J}'_F)\cong
\oplus_{v\in V(L)}H^2(G_w,F^*_w),$$
o pour toute place $v$ de $L$  , $w$ est une place de $F$ au-dessus de $v$.
\medskip

La thorie du corps de classes fournit des isomorphismes canoniques 
$$\widehat H^2(G,C_F)\cong 1/\# G\ {\bf Z}/\ {\bf Z}$$ et 
$$\widehat H^2(G_w,F_w^*)\cong 1/\# G_w\ {\bf Z}/\ {\bf Z}.$$ 

Compte tenu de ces isomorphismes, l'application $s_2$ est donne par 
$$s_2((z_v)_{v\in V(F/L)})=\sum_{v\in V(F/L)}z_v.$$ 

Dans cette galit, l'expression $\sum_{v\in V(F/L)}z_v$ a un sens car l'ordre de $G_w$ divise l'ordre de $G$.

D'autre part, la suite exacte courte de valuation conduit  l'injectivit de l'homomorphisme 
$$\widehat H^2(G_w,A_w^*)\to \widehat H^2(G_w,F^*_w),$$ ce qui montre que $j_2$ est injective.
Par consquent, l'homomorphisme $\varphi_2$ est donne par 
$$\varphi_2((z_v)_{v\in {\cal R}(F/L)})=
\sum_{v\in {\cal R}(F/L)}z_v.$$ 

En particulier, $\ker \varphi_2\not=\widehat H^2(G,{\cal U}_F)$.

Puisque $G$ est cyclique, il existe 
$\gamma\in H^2(G,{\bf Z})$ tel que les homomorphismes verticaux $\gamma\cup-$ du diagramme  ci-dessous sont des isomorphismes.
$$\xymatrix{
\widehat H^0(G,{\cal U}_F)\ar[r]^{\varphi_0}\ar[d]_{\gamma\cup-}&
\widehat H^0(G,C_F)\ar[d]^{\gamma\cup-}\cr
\widehat H^2(G,{\cal U}_F)\ar[r]_{\varphi_2}&
\widehat H^2(G,C_F)\cr
}$$
Par naturalit des cup produits, ce diagramme est commutatif.
 
De ce diagramme et de $\ker\varphi_2\not=\widehat H^2(G,{\cal U}_F)$, on dduit 
$\ker\varphi_0\not=\widehat H^0(G,{\cal U}_F)$. Mais $\widehat H^0(G,C_F)$ est isomorphe   $G$ qui est simple, car d'ordre premier. Par consquent,  $\varphi_0$ est surjectif et on a la suite exacte
 $$\xymatrix{1\ar[r]&
 \hbox{im\ }h_0\ \cap\ \hbox{im\ }j_0\ar[r]&
 \hbox{im\ }j_0\ar[r]&
 G\ar[r]&1.
 }$$ 
Par ailleurs, d'aprs 1)\  $\hbox{im\ }j_0$ est isomorphe  $\widehat H^0(G,{\cal U}_F)$. On a vu  que ce dernier groupe
est isomorphe  
$\oplus_{v\in{\cal R}(F/L)}\widehat H^0(G_w,A^*_w)$. D'aprs [Se], corollaire 7 et sa remarque p. 95, $\widehat H^0(G_w,A^*_w)$ est cyclique d'ordre $\ell$ et donc 
$\hbox{im\ }j_0\cong {\bf C}_\ell^t$ est un ${\bf F}_\ell$-espace vectoriel de dimension $t$. La suite exacte ci-dessus implique donc l'galit suivante :   
$$\hbox{dim\ }_{{\bf F}_\ell}\left(\hbox{im\ }h_0\ \cap\ \hbox{im\ }j_0\right)\  =t-1.$$
 
\bigskip

{\bf 6. Quelques applications.}
\medskip

Nous supposons  prsent et jusqu' la fin de cet article que  $F/L$ est une extension cyclique de corps de nombres, d'anneaux d'entiers $A={\cal O}_F$, 
 $B={\cal O}_L$. Rappelons que $I_GCl(A)$ dsigne le sous-module de $Cl(A)$ engendr par les lments de la forme $(1-g)a$, $a\in Cl(A)$ et que $${}_NCl(A)=\ker N: Cl(A)\to Cl(B).$$ 
 
 \medskip

{\bf Lemme 6.1.}
{\sl a) Avec les notations de 5.5, dans le cas d'une extension cyclique, on a un isomorphisme de groupes 
$$\hbox{coker\ }\alpha\ \cong\  {}_NCl(A)/I_GCl(A).$$

b) Pour une extension cyclique $F/{\bf Q}$, 
$\hbox{coker\ }\alpha$ est  isomorphe au groupe des co-invariants du groupe des classes.}

\medskip

{\sl Dmonstration}. 
a) D\'esignons par $g$ un g\'en\'erateur du groupe de Galois $G$.
Soit $z$ un \'el\'ement de ${\bf J}_F$. Dsignons par 
$I_z$ l'idal fractionnaire $\prod_w\mathfrak p_w^{v(z_w)}$.
Le morphisme de groupes 
$$\overline \alpha:{\bf J}_F\to Cl(A)$$
dfini par 
$\overline\alpha(z)=[I_z]$
est surjectif et satisfait la relation 
$\overline\alpha(gz)=g\overline\alpha(z)$ pour tout $g$ de $G$ et tout $z$ de ${\bf J}_F$.

La restriction $\alpha$ de $\overline \alpha$  
${}_N{\bf J}_F$ est  valeurs dans 
${}_NCl(A)$.
Soit 
$z\in {}_N{\bf J}_F$. 
D'apr\`es le th\'eor\`eme 90 de Hilbert sous la forme 
$H^1(G,{\bf J}_F)=1$, il existe un \'el\'ement $u$ de ${\bf J}_F$ tel que $z=(1-g)u$. On en dduit $\alpha(z)=(1-g)[I_u]$ donc $\alpha(z)\in I_GCl(A)$. Ceci montre l'inclusion $\hbox{im\ }\alpha\subset I_GCl(A)$. 

Soit $(1-g)[I]$ un lment de $I_GCl(A)$. Par surjectivit de $\overline\alpha$, il existe $u\in {\bf J}_F$ tel que $[I]=[I_u]$. Posons $z=(1-g)u$. On a $z\in{}_NCl(A)$ et $\alpha(z)=(1-g)[I]$. Ceci montre l'inclusion 
$I_GCl(A)\subset \hbox{im\ }\alpha$.

b) Si $L={\bf Q}$, ${}_NCl(A)=Cl(A)$ et 
$\hbox{coker\ }(\alpha)=Cl(A)/I_GCl(A)=Cl(A)_G$.

\medskip
 
{\bf Proposition 6.2.}
{\sl Si le groupe de Galois $G$ est cyclique, les homomorphismes naturels 
$$h:L^*/B^*N(F^*)\to {\bf J}_L/B^*N({\bf J}_F)$$
et
$$j:{\cal U}_L/B^*N({\cal U}_F)\to {\bf J}_L/B^*N({\bf J}_F)$$
 
 sont injectifs.}

\medskip

{\sl Dmonstration}. Consid\'erons le diagramme suivant dans lequel les lignes sont  exactes.

$$\xymatrix{
1\ar[r]&
B^*/B^*\cap N({\cal U}_F)\ar[r]\ar[d]_{\tilde j_0}&
\widehat H^0(G,{\cal U}_F)\ar[r]^{p_{\cal U}}
\ar[d]_{j_0}&
{\cal U}_L/B^*N({\cal U}_F)\ar[r]
\ar[d]^{j}&1\cr
1\ar[r]&
B^*/B^*\cap N({\bf J}_F)\ar[r]&
\widehat H^0(G,{\bf J}_F)\ar[r]^p&
{\bf J}_L/B^*N({\bf J}_F)\ar[r]&1\cr
1\ar[r]&
B^*/B^*\cap N(F^*)\ar[r]\ar[u]^{\tilde h_0}&
\widehat H^0(G,F^*)\ar[r]^{p_F}\ar[u]^{h_0}&
L^*/B^*N(F^*)\ar[r]\ar[u]_{h}&1\cr
}$$

Dans ce diagramme, les homomorphismes $\tilde h_0$ et $\tilde j_0$ sont  surjectifs. L'homomorphisme $j_0$ est  injectif et, gr\^ace au principe de Hasse, il en est de m\^eme pour  l'homomorphisme $h_0$ .  En appliquant le lemme du serpent aux deux premi\`eres lignes, on dmontre  la fois l'injectivit de 
$\tilde j_0$ et celle de  $j$. On raisonne de m\^eme pour $\tilde h_0$ et $h$. 
\medskip

Rappelons sans dmonstration un rsultat banal.
\medskip

{\bf Lemme 6.3.}
\smallskip

{\sl Soient $M_1$ et $M_2$ deux sous-groupes d'un groupe ablien $M$ et soit $N$ un sous-groupe de $M_1\cap M_2$. Alors on a un isomorphisme de groupes 
$M_1\cap M_2/N\cong M_1/N\cap M_2/N$.}

\medskip

{\bf Lemme 6.4.}
{\sl Avec les notations de 5.5 et 5.7, dans le cas d'une extension cyclique, on a des isomorphismes de groupes 
$$\ker\mu\ \cong 
\ \hbox{im\ }h_0\ \cap \ \hbox{im\ }j_0$$ 
et
$$\ker\mu/\hbox{im\ }\beta\ \ \cong \ 
\hbox{im\ }h\ \cap\ \hbox{im\ }j.$$}
\medskip

{\sl Dmonstration}.
En notation multiplicative, l'application 
$$\mu:\widehat H^0(G,F^*)\times\widehat H^0(G,{\cal U}_F)\to
\widehat H^0(G,{\bf J}_F)$$ 
de la suite exacte de Mayer-Vietoris est donn\'ee par $$\mu(x,y)=h_0(x)j_0^{-1}(y).$$
 Son noyau $\ker\mu$ est donc isomorphe au produit fibr $X\times_Z Y$ avec 
$X=\widehat H^0(G,F^*)$, 
$Y=\widehat H^0(G,{\cal U}_F)$,
$Z=\widehat H^0(G,{\bf J}_F)$. L'extension $F/L$ tant cyclique, on a vu que les homomorphismes $h_0$ et $j_0$ sont injectifs. Par consquent, ce produit fibr n'est autre que  
$\hbox{im\ }h_0\cap\hbox{im\ }j_0.$ Ceci montre 
$$\ker\mu\cong  \hbox{im\ }h_0\ \cap\ \hbox{im\ }j_0.$$
\smallskip

Considrons le diagramme suivant

$$\xymatrix{
\widehat H^0(G,A^*)\ar[rd]_{\beta}\ar[rrd]^{\beta_{\cal U}}\ar[ddr]_{\beta_F}
&&&\hbox{coker\ }\beta_{\cal U}\ar[ddd]^{j}\cr
&\hbox{ker\ }\mu\ar[r]\ar[d]&\widehat H^0(G,{\cal U}_F)\ar[ur]\ar[d]^{j_0}&\cr
&\widehat H^0(G, F^*)\ar[r]_{h_0}\ar[dl]&\widehat H^0(G,{\bf J}_F)\ar[dr]\cr
\hbox{coker\ }\beta_F\ar[rrr]_{h}&&&\hbox{coker\ }(h_0\circ\beta_F)
}$$

On a 
$\ker\mu\cong \hbox{im\ }h_0\cap \hbox{im\ }j_0,$ d'o
$$\ker\mu/\hbox{im\ }\beta\cong
(\hbox{im\ }h_0\cap \hbox{im\ }j_0)/
\hbox{im\ }(h_0\beta_F)\cong$$

$$\cong \hbox{im\ }h_0/\hbox{im\ }(h_0\beta_F)\ \cap\ \hbox{im\ }j_0/\hbox{im\ }(h_0\beta_{\cal U})\ =\ 
\hbox{im\ }h\cap \hbox{im\ }j.$$

\medskip

{\bf Remarque 6.5.} On a $\hbox{coker\ }\beta_F=L^*/B^*N(F^*)$,
$\hbox{coker\ }\beta_{\cal U}={\cal U}_L/B^*N({\cal U}_F)$
et 
$\hbox{coker\ }h_0\circ\beta={\bf J}_L/B^*N({\bf J}_F).$  

Les homomorphismes $h$ et $j$ tant injectifs, nous nous autoriserons  crire $\hbox{im\ } h\cap \hbox{im\ }  j$ sous la forme
$$\hbox{im\ }  h\cap \hbox{im\ }  j\ \cong\  L^*/B^*N(F^*)\ \cap\ {\cal U}_L/B^*N({\cal U}_F),$$
vu comme sous-groupe de ${\bf J}_L/B^*N({\bf J}_F).$

\medskip

{\bf Thorme 6.6.}
{\sl Soit $F/L$ est une extension cyclique de corps de nombres d'anneaux d'entiers $A={\cal O}_F$, 
 $B={\cal O}_L$. On dsigne par
 $I_GCl(A)$  le sous-module de $Cl(A)$ engendr par les lments de la forme $(1-g)a$, $a\in Cl(A)$ et on pose ${}_NCl(A)=\ker N: Cl(A)\to Cl(B).$  On a alors un isomorphisme de groupes 
$${}_NCl(A)/I_GCl(A)\ \cong\  
\hbox{im\ }h\ \cap\ \hbox{im\ }j.\leqno (6.6.)$$

En particulier, pour une extension $F/{\bf Q}$ cyclique, si on dsigne  par $Cl(A)_G$ le groupe des co-invariants du groupe des classes de $A$ et par 
$\hat {\bf Z}^*$ le produit $\prod_{p} {\bf Z}_p^*$,  on a un isomorphisme
  de groupes
$$Cl(A)_G\ \cong\  
{\bf Q}^*/{\bf Z}^*N(F^*)\ \cap\ \    \widehat {\bf Z}^*/{\bf Z}^*N({\cal U}_F),$$
 ce dernier tant vu comme  sous-groupe de 
 ${\bf J}_{\bf Q}/{\bf Z}^*N({\bf J}_F).$}

\medskip

La dmonstration rsulte immdiatement de 5.5, 6.1 et 6.4.

\medskip

{\bf Corollaire 6.7.}
{\sl Sous les hypothses 6.6 et avec les notations de 5.7, on a une suite exacte 
 
 $$\xymatrix{
 1\ar[r]& 
 B^*/(B^*\cap N(F^*))\ar[r]^-{}&
\hbox{im\ }h_0\ \cap\ \hbox{im\ }j_0\ar[r]&
  {}_NCl(A)/I_GCl(A)  \ar[r]& 1.
 }$$}
  
 {\sl Dmonstration}. Dans le diagramme de la  dmonstration 6.2, on v\'erifie par chasse dans le diagramme que la restriction de $p$ \`a 
 $\hbox{im\ }h_0\ \cap\  \hbox{im\ }j_0$ a pour image 
 $\hbox{im\ }h\ \cap\  \hbox{im\ }j$ et pour noyau 
  $\ker p=B^*/(B^*\cap N(F^*))$. On a donc une suite exacte courte 
$$\xymatrix{
 1\ar[r]& 
 B^*/(B^*\cap N(F^*))\ar[r]^-{}&
 \hbox{im\ }h_0\ \cap\ \hbox{im\ }j_0\ar[r]&
  \hbox{im\ }h\ \cap\ \hbox{im\ }j\ar[r]& 0.
 }$$
 
 On conclut gr\^ace  6.6.
\medskip
  
  De 6.7 et 5.7.4, nous dduisons immdiatement 
  
  \medskip
  
{\bf Thorme 6.8.}
{\sl Soit $F/L$ une extension cyclique de corps de nombres d'anneaux d'entiers $A={\cal O}_F$ et $B={\cal O}_L$. On suppose l'extension ramifie, de degr premier $\ell$ et de groupe de Galois $G$. Dsignons par $t$ le nombre de places finies de $L$ ramifies dans $F$ et par
$I_GCl(A)$ le sous-module de $Cl(A)$ engendr par les lments de la forme $(1-g)a$, $a\in Cl(A)$. Posons  
$${}_NCl(A):=\ker N:Cl(A)\to Cl(B).$$ Alors 
${}_NCl(A)/I_GCl(A)$ est un ${\bf F}_\ell$-espace vectoriel de dimension infrieure ou gale  $t-1$.}
  
 \medskip

{\sl Exemple d'emploi de la suite exacte 6.7.}
 \medskip

 Montrons  comment la suite exacte 6.7 permet de retrouver un rsultat  connu   sur le $\ell$-rang des extensions de degr premier $\ell$, mais de dmonstration dlicate. 
  \medskip

{\bf Th\'eor\`eme 6.9.}
{\sl Soit $F/{\bf Q}$ une extension cyclique de degr\'e premier $\ell$, de groupe de Galois $G={\bf C}_\ell$. Notons $Cl(A)_G$ le groupe des co-invariants du groupe des classes de $F$ sous l'op\'eration de son groupe de Galois.
Soit $t$ le nombre de diviseurs premiers du discriminant du corps $F$.

Alors le $\ell$-rang du groupe des co-invariants du groupe des classes du corps $F$ est \'egal \`a $t-1$ sauf si $\ell=2$, $F$ quadratique r\'eel et s'il existe un nombre premier $p$ congru \`a $3$ modulo $4$ et  divisant le discriminant du corps $F$. 

Dans ce dernier cas,  le $2$-rang du groupe des co-invariants du groupe des classes du corps $F$ est \'egal \`a $t-2$
.}
\medskip

{\sl Dmonstration}. 
Posons   
$V=\hbox{im\ }h_0\ \cap\ \hbox{im\ }j_0$
et
$H={\bf Z}^*/{\bf Z}^*\cap N(F^*)$.
D'aprs 5.7.4, $V$ est un ${\bf F}_\ell$-espace vectoriel  de dimension $t-1$.
D'apr\`es 6.7, on a une suite exacte courte de ${\bf F}_\ell$-espaces vectoriels
$$\xymatrix{
1\ar[r]& H\ar[r]^-{}&
V\ar[r]&
Cl(A)_G\ar[r]& 1.}$$

Si $\ell>2$ ou si $l=2$ et $-1\in N(F^*)$, $H=\{1\}$,  $V\cong Cl(A)_G$ et $$\dim_{{\bf F}_\ell} Cl(A)_G=\dim_{{\bf F}_\ell} V=t-1.$$

Rappelons que pour  $\ell=2$, on a    $-1\not\in N(F^*)$ si et seulement   $F$ est r\'eel et s'il existe un premier congru  3 modulo 4 ramifi dans $F$.

Pour $\ell=2$ et $-1\not\in N(F^*)$,  
on a  $H={\bf Z}^*$
donc $\dim_{{\bf F}_\ell} Cl(A)_G=\dim_{{\bf F}_\ell} V-1=t-2$.
\medskip

{\bf Remarque 6.10.} Le thorme 6.9 a pour origine  une formule des genres. On pourra trouver des dmonstrations classiques dans  S. Lang [L], Chap. XIII, lemme 4.1 ou dans [H], Theorem 4, p. VII.12. 
Pour $\ell=2$, on retrouve un rsultat essentiellement d\^u  Gauss et qui s'nonce sous la forme suivante : soit $\Delta$ un discriminant fondamental de corps quadratique et soit 
$t=t(\Delta)$ le nombre de facteurs premiers de $\Delta$. Le $2$-rang du groupe des classes de tout corps quadratique de discriminant $\Delta$ est gal  $t-1$, sauf si $\Delta>0$ et si $\Delta$ possde un facteur premier congru  $3$ modulo $4$. Dans ce cas exceptionnel, le $2$-rang du groupe des classes est gal  $t-2$.
\bigskip
\bigskip

{\bf Bibliographie.}

[B], H. Bass,   {\sl  Algebraic K-theory}, Benjamin, 1968.

[B-M-S], H. Bass, J. Milnor \& J.-P. Serre, Solution of the congruence subgroup problem for $SL_n$ ($n\geq 3$) and $Sp_{2n}$ ($n\geq 2$), {\sl Publ. Math. Inst. Hautes \'Et. Sci.}, {\bf 33}, 1967, 59-137.

[G], S.M. Gersten, {\sl Higher algebraic $K$-theory}, Lecture Notes in Math.,  341, Springer, 1973, 3-41.

 
[Gr], D. Grayson,  Higher algebraic $K$-theory. II (after Daniel Quillen).  {\sl Algebraic $K$-theory (Proc. Conf., Northwestern Univ., Evanston, Ill., 1976)},  pp. 217--240. Lecture Notes in Math., Vol. 551, Springer, Berlin, 1976.

[H], C.S. Herz, {\sl Seminar on Complex Multiplication}, Expos\'e VII,  Lecture Notes in Math.,  21, Springer, 1966.

[K1], M. Karoubi, Foncteurs drivs et $K$-thorie, {\sl Lecture Notes in Math.}, 136, 1970, 107-186.

[K2], M. Karoubi, Localisation de formes quadratiques 1, {\sl Ann. Sci. \'E.N.S.}, {\bf 7}, 1974, 359-404.

[K3], M. Karoubi, Le th\'eor\`eme fondamental de la $K$-th\'eorie hermitienne, {\sl Ann. Math.}, {\bf 112}, 1980, 259-282.

[K-L], M. Karoubi \& T. Lambre, Quelques classes caract\'eristiques en th\'eorie des nombres, {\sl J. reine angew. Math.}, {\bf 543}, 2002, 169-186.

[L], S. Lang,  {\sl Cyclotomic fields}, II, Graduate Texts in Mathematics 69, Springer, 1979.

[M], J. Milnor, {\sl Introduction to algebraic $K$-theory}, Annals of Mathematics Studies, No. 72. Princeton University Press, Princeton, N.J.; University of Tokyo Press, Tokyo, 1971.

[N], J. Neukirch, {\sl Algebraic number theory}, Springer, 1999.

[N-S-W], J. Neukirch, A. Schmidt, K. Wingberg, {\sl Cohomology of number fields}, Springer, 2000.

[NQD] T. Nguyen Quang Do, {\sl Quelques suites exactes en thorie des genres}. {\sl Publ. Math. Univ. Franche-Comt Besan\c con}, 2006, 103-115.

[Q], D. Quillen, {\sl Higher algebraic $K$-theory}, Lecture Notes in Math., 341, 1973, 85-147.

[Se], J.-P. Serre, {\sl Corps Locaux}, Hermann, 1968.

[Sou], C. Soul, $K$-thorie des anneaux d'entiers de corps de nombres et cohomologie tale, {\sl Invent. Math.}, {\bf 55}, 1979,  251-295.

[Wa], J.B. Wagoner, Delooping classifying spaces in algebraic $K$-theory.  {\sl Topology},  {\bf 11},  1972, 349--370. 

[W1], C. Weibel, K-theory and analytic isomorphism, {\sl Inv. Math.}, {\bf 61}, 1980, 177-197.

[W2], {\sl An introduction to algebraic K-theory (a graduate textbook in progress)}, 

http://math.rutgers.edu/~weibel/Kbook.html

\end{document}